\documentclass[a4paper,11pt]{amsart}
\usepackage{amssymb,amsmath}
\usepackage{amsopn}
\usepackage{mathrsfs}
\usepackage{color}
\usepackage{mathtools}
\usepackage{wasysym}
\usepackage{vmargin}
\usepackage{pdfsync}
\usepackage{vmargin}
\usepackage{color}
\usepackage{amscd}
\usepackage{verbatim}
\newtheorem{theorem}{Theorem}[section]
\newtheorem{definition}[theorem]{Definition}
\newtheorem{lemma}[theorem]{Lemma}
\newtheorem{proposition}[theorem]{Proposition}

\numberwithin{equation}{section}

\newcommand{\CC}{\mathbb{C}}

\newcommand{\rr}{\mathbb{R}}
\newcommand{\eps}{\varepsilon}
\newcommand{\nn}{\mathbb{N}}
\newcommand{\cc}{\mathbb{C}}
\def\un{{\mathrm{1~\hspace{-1.4ex}l}}}

\def\N{\mathbb N}

\def\R{\mathbb R}

\def\norm#1{\Vert#1\Vert}
\def\val#1{\vert#1\vert}

\def\l2{L^2(\R^{n})}
\def\L2{L^2(\R^{2n})}
\def\supp{\operatorname{supp}}

\def\mat22#1#2#3#4{\begin{pmatrix}#1&#2\\ #3&#4\end{pmatrix}}

\author{Karine \textsc{Beauchard}, Philippe \textsc{Jaming} \&  Karel \textsc{Pravda-Starov}}

\address{\noindent \textsc{Karine Beauchard, Univ Rennes, CNRS, IRMAR - UMR 6625, F-35000 Rennes, France
}}
\email{karine.beauchard@ens-rennes.fr}

\address{\noindent \textsc{Philippe Jaming, Universit\'e de Bordeaux, Institut de Math\'ematiques de Bordeaux, UMR 5251,
Cours de la Lib\'eration, F-33405 Talence cedex, France}}
\email{philippe.jaming@math.u-bordeaux.fr}

\address{\noindent \textsc{Karel Pravda-Starov, Univ Rennes, CNRS, IRMAR - UMR 6625, F-35000 Rennes, France
}}
\email{karel.pravda-starov@univ-rennes1.fr}


\thanks{The first and third authors express their gratefulness to the Centre de Math\'ematiques Henri Lebesgue for the very stimulating scientific environment.}

\begin{document}

\title[Spectral estimates for Hermite functions and null-controllability]{Spectral estimates for finite combinations of Hermite functions and null-controllability of hypoelliptic quadratic equations}

\begin{abstract}
Some recent works have shown that the heat equation posed on the whole Euclidean space is null-controllable in any positive time if and only if the control subset is a thick set. This necessary and sufficient condition for null-controllability is linked to some uncertainty principles, 
as the Logvinenko-Sereda theorem, which give limitations on the simultaneous concentration of a function and its Fourier transform. 
In the present work, we prove new uncertainty principles for finite combinations of Hermite functions. We establish an analogue of the
Logvinenko-Sereda theorem with an explicit control of the constant with respect to the energy level of the Hermite functions as eigenfunctions of the harmonic oscillator for thick control subsets. This spectral inequality allows to derive the null-controllability in any positive time from thick control regions for parabolic equations associated with a general class of hypoelliptic non-selfadjoint quadratic differential operators. More generally, the spectral estimate for finite combinations of Hermite functions is actually shown to  hold for any measurable control subset of positive Lebesgue measure, and some quantitative estimates of the constant with respect to the energy level are given for two other classes of control subsets including the case of non-empty open control subsets.  
\end{abstract}

\maketitle
 
\section{Introduction}

The classical uncertainty principle was established by Heisenberg. It points out the fundamental problem in quantum mechanics that the position and the momentum of particles cannot be both determined explicitly, but only in a probabilistic sense with an uncertainty. More generally, uncertainty principles are mathematical results that give limitations on the simultaneous concentration of a function and its Fourier transform. 
When using the following normalization for the Fourier transform
\begin{equation}\label{fourier}
\widehat{f}(\xi)=\int_{\rr^n}f(x)e^{-ix \cdot \xi}dx, \quad \xi \in \rr^n,
\end{equation}
the mathematical formulation of the Heisenberg's uncertainty principle can be stated in a directional version as 
\begin{equation}\label{heis}
\inf_{a \in \rr}\Big(\int_{\rr^n}(x_j-a)^2|f(x)|^2dx\Big)\inf_{b \in \rr}\Big(\frac{1}{(2\pi)^n}\int_{\rr^n}(\xi_j-b)^2|\widehat{f}(\xi)|^2d\xi\Big) \geq \frac{1}{4}\|f\|_{L^2(\rr^n)}^4,
\end{equation}
for all $f \in L^2(\rr^n)$ and $1 \leq j \leq n$. It shows that a function and its Fourier transform cannot both be arbitrarily localized. Moreover, the inequality (\ref{heis}) is an equality if and only if $f$ is of the form 
$$f(x)=g(x_1,...,x_{j-1},x_{j+1},...,x_n)e^{-ibx_j}e^{-\alpha (x_j-a)^2},$$
where $g$ is a function in $L^2(\rr^{n-1})$, $\alpha>0$, and $a$ and $b$ are real constants for which the two infima in (\ref{heis}) are achieved. There are various uncertainty principles of different nature. We refer in particular the reader to the survey article by Folland and Sitaram~\cite{folland}, and the book of Havin and J\"oricke~\cite{havin} for detailed presentations and references for these topics.


Another formulation of uncertainty principles is that a non-zero function and its Fourier transform cannot both have small supports.
For instance, a non-zero $L^2(\rr)$-function whose Fourier transform is compactly supported must be an analytic function with a discrete zero set 
and therefore a full support. This leads to the notion
of weak annihilating pairs as well as the corresponding quantitative notion of strong annihilating ones:

\medskip

\begin{definition} [Annihilating pairs]
Let $S,\Sigma$ be two measurable subsets of $\rr^n$. 
\begin{itemize}
\item[-] The pair $(S,\Sigma)$ is said to be a \emph{weak annihilating pair} if the only function $f\in L^2(\rr^n)$
with $\supp f\subset S$ and $\supp\widehat{f} \subset \Sigma$ is zero $f=0$.
\item[-] The pair $(S,\Sigma)$ is said to be a \emph{strong annihilating pair} if there exists a positive constant $C=C(S,\Sigma)>0$ such that 
for all $f \in L^2(\rr^n)$,
\begin{equation}\label{strongly}
\int_{\rr^n}|f(x)|^2dx \leq C\Big(\int_{\rr^n \setminus S}|f(x)|^2dx  
+\int_{\rr^n \setminus \Sigma}|\widehat{f}(\xi)|^2d\xi\Big).
\end{equation}
\end{itemize}
\end{definition}

\medskip

It can be readily checked that a pair $(S,\Sigma)$ is a strong annihilating pair if and only if there exists a positive constant 
$D=D(S,\Sigma)>0$ such that for all $f \in L^2(\rr^n)$ with $\supp\widehat{f} \subset \Sigma$,
\begin{equation}\label{strongly2}
\|f\|_{L^2(\rr^n)} \leq D\|f\|_{L^2(\rr^n \setminus S)}.
\end{equation}
As already mentioned above in the one-dimensional setting, the pair $(S,\Sigma)$ is a weak annihilating one if $S$ and $\Sigma$ are compact sets. More generally, Benedicks has 
shown in~\cite{benedicks} that $(S,\Sigma)$ is a weak annihilating pair if $S$ and $\Sigma$ are sets of finite Lebesgue measure 
$|S|, |\Sigma| <+\infty$. Under this assumption, the result of Amrein-Berthier~\cite{amrein} actually shows that the pair 
$(S,\Sigma)$ is a strong annihilating one. The estimate $C(S,\Sigma)\leq \kappa e^{\kappa |S||\Sigma|}$ (which is sharp up to numerical 
constant $\kappa>0$) has been established by Nazarov~\cite{Na} in dimension $n=1$. This result was extended to the multi-dimensional case by the second author~\cite{Ja}, with the quantitative estimate 
$C(S,\Sigma)\leq \kappa e^{\kappa (|S||\Sigma|)^{1/n}}$ 
holding if in addition one of the two subsets of finite Lebesgue measure $S$ or $\Sigma$ is convex.

An exhaustive description of all strong annihilating pairs seems for now totally out of reach. We refer the reader for instance to the works~\cite{AO,BD1,BD2,De,DJ,SVW} for a large variety of results and techniques available as well as for examples of weak annihilating pairs that are not strong annihilating ones.
However, there is a complete description of all the support sets $S$ forming a strong annihilating pair with any bounded spectral set $\Sigma$. This description is given by the Logvinenko-Sereda theorem \cite{Logvinenko_Sereda}:

\medskip

\begin{theorem}[Logvinenko-Sereda]\label{Sere}
Let $S,\Sigma\subset\rr^n$ be measurable subsets with $\Sigma$ bounded. Denoting $\tilde S=\rr^n\setminus S$,
the following assertions are equivalent:
\begin{itemize}
\item[-] The pair $(S,\Sigma)$ is a strong annihilating pair
\item[-] The subset $\tilde{S}$ is thick, that is, there exists a cube $K \subset \rr^n$ with sides parallel to coordinate axes and a positive constant $0<\gamma \leq 1$ such that 
$$\forall x \in \rr^n, \quad |(K+x) \cap \tilde{S}| \geq \gamma|K|>0,$$
where $|A|$ denotes the Lebesgue measure of the measurable set $A$.
\end{itemize}
\end{theorem}

\medskip

It is interesting to observe that if $(S,\Sigma)$ is a strong annihilating pair for some bounded subset $\Sigma$, then $S$ makes up a strong annihilating pair with every bounded subset $\Sigma$, but the above constants $C(S,\Sigma)>0$ and $D(S,\Sigma)>0$ do depend on $\Sigma$.
In order to use this result in the control theory of partial differential equations, it is essential to understand how the positive constant $D(S,\Sigma)>0$ depends on the Lebesgue measure of the bounded set $\Sigma$. This question was answered by Kovrijkine~\cite{thesis,Kovrijkine} who established the following quantitative estimates:

\begin{theorem}[Kovrijkine] There exists a universal positive constant $C_n>0$ depending only on the dimension $n\geq 1$ such that
if $\tilde{S}$ is a $\gamma$-thick set at scale $L>0$, that is, for all $x \in \rr^n$,
\begin{equation}\label{thick1v}
|\tilde{S} \cap (x+[0,L]^n)| \geq \gamma L^n,
\end{equation}
with $0<\gamma \leq 1$,
then we have for all $R>0$ and $f \in L^2(\rr^n)$ with 
$\supp\widehat{f} \subset [-R,R]^n$,
\begin{equation}\label{kov}
\|f\|_{L^2(\rr^n)} \leq \Big(\frac{C_n}{\gamma}\Big)^{C_n(1+LR)}\|f\|_{L^2(\tilde{S})}.
\end{equation}
\end{theorem}

Thanks to this explicit dependence of the constant with respect to the parameter $R>0$ in (\ref{kov}), Egidi  and Veseli\'c~\cite{veselic}, and Wang, Wang, Zhang and Zhang~\cite{Wang} have independently established that the heat equation 
\begin{equation}\label{heat}
\left\lbrace \begin{array}{ll}
(\partial_t -\Delta_x)f(t,x)=u(t,x)\un_{\omega}(x), \quad &  x \in \mathbb{R}^n,\ t>0, \\
f|_{t=0}=f_0 \in L^2(\rr^n),                                       &  
\end{array}\right.
\end{equation}
is null-controllable in any positive time $T>0$ from a measurable control subset $\omega \subset \rr^n$ if and only if this subset $\omega$ is thick in $\rr^n$. The notion of null-controllability is defined as follows:

\medskip

\begin{definition} [Null-controllability] Let $P$ be a closed operator on $L^2(\rr^n)$ which is the infinitesimal generator of a strongly continuous semigroup $(e^{-tP})_{t \geq 0}$ on $L^2(\rr^n)$, $T>0$ and $\omega$ be a measurable subset of $\mathbb{R}^n$. 
The equation 
\begin{equation}\label{syst_general}
\left\lbrace \begin{array}{ll}
(\partial_t + P)f(t,x)=u(t,x)\un_{\omega}(x), \quad &  x \in \mathbb{R}^n,\ t>0, \\
f|_{t=0}=f_0 \in L^2(\rr^n),                                       &  
\end{array}\right.
\end{equation}
is said to be {\em null-controllable from the set $\omega$ in time} $T>0$ if, for any initial datum $f_0 \in L^{2}(\mathbb{R}^n)$, there exists $u \in L^2((0,T)\times\mathbb{R}^n)$, supported in $(0,T)\times\omega$, such that the mild (or semigroup) solution of \eqref{syst_general} satisfies $f|_{t=T}=0$.
\end{definition}

\medskip

By the Hilbert Uniqueness Method, see \cite[Theorem~2.44]{coron_book} or \cite{JLL_book}, the null controllability of the equation \eqref{syst_general} is equivalent to the (final state) observability of the adjoint system 
\begin{equation} \label{adj_general}
\left\lbrace \begin{array}{ll}
(\partial_t + P^*)g(t,x)=0, \quad & x \in \mathbb{R}^n,\ t>0, \\
g|_{t=0}=g_0 \in L^2(\rr^n),
\end{array}\right.
\end{equation}
where $P^*$ stands for the $L^2(\rr^n)$-adjoint of $P$.
This notion of observability is defined as follows:

\medskip

\begin{definition} [Observability] Let $T>0$ and $\omega$ be a measurable subset of $\mathbb{R}^n$. 
Equation \eqref{adj_general} is said to be {\em observable from the set $\omega$ in time} $T>0$ if there exists a positive constant $C_T>0$ such that,
for any initial datum $g_0 \in L^{2}(\mathbb{R}^n)$, the mild (or semigroup) solution of \eqref{adj_general} satisfies
\begin{equation}\label{eq:observability}
\int\limits_{\mathbb{R}^n} |g(T,x)|^{2} dx  \leq C_T \int\limits_{0}^{T} \Big(\int\limits_{\omega} |g(t,x)|^{2} dx\Big) dt.
\end{equation}
\end{definition}

\medskip

Following~\cite{veselic} or~\cite{Wang}, the necessity of the thickness property of the control subset for the null-controllability in any positive time is a consequence of a quasimodes construction; whereas the sufficiency is derived in~\cite{veselic} from an abstract observability result obtained by an adapted Lebeau-Robbiano method~\cite{robbia} and established by the first and third authors with some contributions of Luc Miller\footnote{Universit\'e Paris-Ouest, Nanterre La D\'efense, UFR SEGMI, B\^atiment G, 200 Av. de la R\'epublique, 92001 Nanterre Cedex, France (luc.miller@math.cnrs.fr)}:

\medskip

\begin{theorem} \label{Meta_thm_AdaptedLRmethod}
\cite[Theorem~2.1]{KK1}. Let $\Omega$ be an open subset of $\mathbb{R}^n$,
 $\omega$ be a measurable subset of $\Omega$,
 $(\pi_k)_{k \in \mathbb{N}^*}$ be a family of orthogonal projections defined on $L^2(\Omega)$,
 $(e^{-tA})_{t \geq 0}$ be a strongly continuous contraction semigroup on $L^2(\Omega)$;
 $c_1, c_2, a, b, t_0, m >0$ be positive constants with $a<b$.
If the following spectral inequality
\begin{equation} \label{Meta_thm_IS}
\forall g \in L^2(\Omega), \forall k \geq 1, \quad \|\pi_k g \|_{L^2(\Omega)} \leq e^{c_1 k^a} \|\pi_k g \|_{L^2(\omega)},
\end{equation}
and the following dissipation estimate 
\begin{equation} \label{Meta_thm_dissip}
\forall g \in L^2(\Omega), \forall k \geq 1, \forall 0<t<t_0, \quad \| (1-\pi_k)(e^{-tA} g)\|_{L^2(\Omega)} \leq \frac{1}{c_2} e^{-c_2 t^m k^b} \|g\|_{L^2(\Omega)},
\end{equation}
hold, then there exists a positive constant $C>1$ such that the following observability estimate holds
\begin{equation} \label{meta_thm_IO}
\forall T>0, \forall g \in L^2(\Omega), \quad \| e^{-TA} g \|_{L^2(\Omega)}^2 \leq C\exp\Big(\frac{C}{T^{\frac{am}{b-a}}}\Big) \int_0^T \|e^{-tA} g \|_{L^2(\omega)}^2 dt.
\end{equation}
\end{theorem}

\medskip

The proof of the above result is inspired by the works~\cite{mi1,mi2}.
In the statement of~\cite[Theorem~2.1]{KK1}, the subset $\omega$ is supposed to be open in $\Omega$. However, the proof given in~\cite{KK1} works as well when the subset $\omega$ is only assumed to be measurable.
Notice that the assumptions in the above statement do not require that the orthogonal projections $(\pi_k)_{k \geq 1}$ are spectral projections onto the eigenspaces of the infinitesimal generator $A$, which is allowed to be non-selfadjoint.

According to the above result, there are two key ingredients to derive a result of null-controllability, or equivalently a result of observability, while using Theorem~\ref{Meta_thm_AdaptedLRmethod}: a spectral inequality (\ref{Meta_thm_IS}) and a dissipation estimate (\ref{Meta_thm_dissip}). For the heat equation, the orthogonal projections used are the frequency cutoff operators given by the orthogonal projections onto the closed vector subspaces
\begin{equation}\label{cutoff}
E_k=\{f \in L^2(\rr^n) : \textrm{ supp }\widehat{f} \subset [-k,k]^n \}, \end{equation}
for $k \geq 1$.
With this choice, the dissipation estimate readily follows from the explicit formula
\begin{equation}\label{heat1}
\widehat{(e^{t\Delta_x}g)}(t,\xi)=\widehat{g}(\xi)e^{-t|\xi|^2}, \quad t \geq 0, \ \xi \in \rr^n,
\end{equation}     
whereas the spectral inequality is given by the sharp formulation of the Logvinenko-Sereda theorem (\ref{kov}). Notice that the power $1$ for the parameter $R$ in (\ref{kov}) and the power $2$ for the term $|\xi|$ in (\ref{heat1}) account for the fact that Theorem~\ref{Meta_thm_AdaptedLRmethod} can be applied with the parameters $a=1$, $b=2$ that satisfy the required condition $0<a<b$. It is therefore essential that the power of the parameter $R$ in the exponent of the estimate (\ref{kov}) is strictly less than $2$. As there is still a gap between the cost of the localization ($a=1$) given by the spectral inequality and its compensation by the dissipation estimate ($b=2$), it is interesting to notice that we could have expected that the null-controllability of the heat equation could have held under weaker assumptions than the thickness property on the control subset, by allowing some higher costs for localization with some parameters $1<a<2$, but the Logvinenko-Sereda theorem actually shows that this is not the case. Indeed, if the spectral inequality holds with a parameter $1<a<2$, the pair $(\rr^n \setminus \omega, [-k,k]^n)$ has to be a strong annihilating one and $\omega$ has to be thick according to Theorem~\ref{Sere}.

Theorem~\ref{Meta_thm_AdaptedLRmethod} does not only apply with the use of frequency cutoff projections and a dissipation estimate induced by some Gevrey-type regularizing effects\footnote{Following~\cite{morixu}, the Gevrey-type spaces $\mathcal{A}^s(\rr^n)$, with $s>0$, are defined as the spaces of smooth functions $f \in C^{\infty}(\rr^n)$ satisfying
$$\exists C>1, \forall \alpha \in \nn^n, \quad \|\partial_x^{\alpha}f\|_{L^2(\rr^n)} \leq C^{1+|\alpha|}(\alpha!)^s.$$
Thanks to Sobolev embeddings, the $L^2$-norm can be replaced by the $L^{\infty}$-one in the above estimates.}. 
Other regularities than the Gevrey one can be taken into account. In the previous work~\cite{KK1} by the first and third authors,  Theorem~\ref{Meta_thm_AdaptedLRmethod} is used for a general class of accretive hypoelliptic quadratic operators $q^w$ generating some strongly continuous contraction semigroups $(e^{-tq^w})_{t \geq 0}$ enjoying some Gelfand-Shilov regularizing effects. The definition and standard properties related to Gelfand-Shilov regularity\footnote{As explained in Section~\ref{quadratic1}, this notion of regularity plays a key role to obtain the results of null-controllability in Theorem~\ref{th2} as an application of abstract results established in~\cite{KK1}, but is not used to derive the uncertainty principles in Theorem~\ref{th1}.} are recalled in Appendix (Section~\ref{gelfand}). As recalled in this appendix, the Gelfand-Shilov regularity is characterized by specific exponential decays of the functions and their Fourier transforms; and in the symmetric case, can be read on the exponential decay of the Hermite coefficients of the functions in their expansions in the $L^2(\rr^n)$-Hermite basis $(\Phi_{\alpha})_{\alpha \in \nn^n}$.
Explicit formulas and some reminders of basic facts about Hermite functions are given in Section~\ref{hermite_functions}. The class of hypoelliptic quadratic operators whose description will be given in Section~\ref{quadratic1} enjoys some Gelfand-Shilov regularizing effects ensuring that the following dissipation estimate holds~\cite[Proposition~4.1]{KK1}: 
\begin{multline}\label{eq6}
\exists C_0>1, \exists t_0>0, \forall t > 0, \forall k \geq 0, \forall f \in L^2(\rr^n),\\
\|(1-\pi_k)(e^{-tq^w}f)\|_{L^2(\rr^n)} \leq C_0e^{- \delta(t)k} \|f\|_{L^2(\rr^n)},
\end{multline} 
with 
\begin{equation}\label{eq5}
\delta(t)=\frac{\inf(t,t_0)^{2k_0+1}}{C_0} > 0, \quad t > 0, \qquad 0 \leq k_0 \leq 2n-1,
\end{equation} 
where  
\begin{equation}\label{estt3}
\mathbb{P}_k g =\sum_{\substack{\alpha \in \N^n,\\ |\alpha|=k}}\langle g,\Phi_{\alpha}\rangle_{L^2(\rr^n)}\Phi_{\alpha}, \quad k \geq 0, 
\end{equation} 
with $\val \alpha=\alpha_{1}+\dots+\alpha_{n}$,
denotes the orthogonal projection onto the $k^{\textrm{th}}$ energy level associated with the harmonic oscillator 
$$\mathcal{H}=-\Delta_x+|x|^2=\sum_{k=0}^{+\infty}(2k+n)\mathbb{P}_k,$$ 
and
\begin{equation}\label{estt4}
\pi_k=\sum_{j=0}^k\mathbb{P}_j, \quad k \geq 0,
\end{equation}
denotes the orthogonal projection onto the $(k+1)^{\textrm{th}}$ first energy levels. The above integer $0 \leq k_0 \leq 2n-1$ whose definition is given in (\ref{h1bis2}) is a structural parameter intrinsically associated to the Weyl symbol of the quadratic operator. 
In order to apply Theorem~\ref{Meta_thm_AdaptedLRmethod}, we need a spectral inequality for finite combinations of Hermite functions of the type 
\begin{equation}\label{eq7}
\exists C>1, \forall k \geq 0, \forall f \in L^2(\mathbb{R}^n), \quad \|\pi_k f\|_{L^2(\rr^n)} \leq Ce^{Ck^a} \|\pi_k f\|_{L^2(\omega)},
\end{equation}
with $a<1$, where $\pi_k$ is the orthogonal projection (\ref{estt4}). In~\cite[Proposition~4.2]{KK1}, such a spectral inequality is established with $a=\frac{1}{2}$ when the control subset $\omega$ is an open subset of $\rr^n$ satisfying the following geometric condition:  
\begin{equation}\label{hyp_omega}
\exists \delta, r >0, \forall y \in \mathbb{R}^n\,, \exists y' \in \omega,\quad 
\left\lbrace \begin{array}{ll}
B(y',r) \subset \omega, \\ 
|y-y'|<\delta,
\end{array}\right.
\end{equation}
where $B(y',r)$ denotes the open Euclidean ball in $\rr^n$ centered at $y'$ with radius $r>0$. It allows to derive the null-controllability of parabolic equations associated to accretive quadratic operators with zero singular spaces in any positive time $T>0$ from any open subset $\omega$ of $\rr^n$ satisfying (\ref{hyp_omega}). The above geometric condition was introduced by Le~Rousseau and Moyano in~\cite{moyano} and was shown to be sufficient for the null-controllability of the Kolmogorov equation in any positive time.

In the present work, we study under which geometric conditions on the control subset~$\omega \subset \rr^n$, the spectral inequality
\begin{equation}\label{eq7l}
\forall k \geq 0, \exists C_k(\omega)>0, \forall f \in L^2(\mathbb{R}^n), \quad \|\pi_k f\|_{L^2(\rr^n)} \leq C_k(\omega) \|\pi_k f\|_{L^2(\omega)},
\end{equation}
holds and how the geometric properties of the set $\omega$ relate to the possible growth of the positive constant $C_k(\omega)>0$ with respect to the energy level when $k \to +\infty$. The main results contained in this article provide some
quantitative upper bounds on the positive constant $C_k(\omega)>0$ with respect to the energy level for three different classes of measurable subsets~: 
\begin{itemize}
\item[-] non-empty open subsets in $\rr^n$,
\item[-] measurable sets in $\rr^n$ verifying the condition
\begin{equation}\label{liminf}
\liminf_{R \to +\infty}\frac{|\omega \cap B(0,R)|}{|B(0,R)|}=\lim_{R \to +\infty}\Big(\inf_{r \geq R}\frac{|\omega \cap B(0,r)|}{|B(0,r)|}\Big)>0,
\end{equation}
where $B(0,R)$ denotes the open Euclidean ball in $\rr^n$ centered at $0$ with radius $R>0$,
\item[-] thick measurable sets in $\rr^n$.  
\end{itemize}
We observe that in the first two classes, the measurable control subsets are allowed to have gaps containing balls with radii tending to infinity, whereas in the last class there must be a bound on such radii. We shall see that the quantitative upper bounds obtained for the two first classes (Theorem~\ref{th1}, estimates $(i)$ and $(ii)$) are not sufficient to obtain any result of null-controllability for the class of hypoelliptic quadratic operators studied in Section~\ref{quadratic1}. 
Regarding the third one, the quantitative upper bound (Theorem~\ref{th1}, estimate $(iii)$) is an analogue of the Logvinenko-Sereda theorem for finite combinations of Hermite functions. As an application of this third result, we extend in Theorem~\ref{th2} the result of null-controllability for parabolic equations associated to accretive quadratic operators with zero singular spaces from any thick set $\omega \subset \rr^n$ in any positive time $T>0$.

\section{Uncertainty principles for finite combinations of Hermite functions}
Let $(\Phi_{\alpha})_{\alpha \in \nn^n}$ be the $n$-dimensional Hermite functions whose definitions are recalled in Section~\ref{hermite_functions} and 
\begin{equation}\label{jk1b}
\mathcal E_{N}=\text{Span}_{\cc}\{\Phi_{\alpha}\}_{\alpha\in \N^n, \val \alpha \leq N},
\end{equation}
be the finite dimensional vector space spanned by all the Hermite functions with $|\alpha| \leq N$.

As the Lebesgue measure of the zero set of a non-zero analytic function on $\cc$ is zero, the $L^2$-norm $\|\cdot\|_{L^2(\omega)}$ on any measurable set $\omega \subset \rr$ of positive measure $|\omega|>0$ defines a norm on the finite dimensional vector space $\mathcal E_{N}=\pi_{N}(L^2(\rr^n))$, with $\pi_N$ the orthogonal projection defined in (\ref{estt4}). As a consequence of the Remez inequality, we check in Appendix (Section~\ref{remez0}) that this result holds true as well in the multi-dimensional case when  $\omega \subset \rr^n$, with $n \geq 1$, is a measurable subset of positive Lebesgue measure $|\omega|>0$. By equivalence of norms in finite dimension, for any measurable set $\omega \subset \rr^n$ of positive Lebesgue measure $|\omega|>0$ and all $N \in \nn$, there exists a positive constant $C_N(\omega)>0$ depending on $\omega$ and $N$ such that the following spectral inequality holds
\begin{equation}\label{spec}
\forall f \in L^2(\rr^n), \quad \|\pi_Nf\|_{L^2(\rr^n)} \leq C_N(\omega)\|\pi_Nf\|_{L^2(\omega)}.
\end{equation} 
We aim at studying how the geometric properties of the set $\omega$ relate to the possible growth of the positive constant $C_N(\omega)>0$ with respect to the energy level.

The main results of the present work are given by the following uncertainty principles for finite combinations of Hermite functions:

\medskip  

\begin{theorem}\label{th1}
The following spectral estimates hold:

\medskip

\noindent
$(i)$ If $\omega$ is a non-empty open subset of $\rr^n$, then there exists a positive constant $C=C(\omega)>1$ such that  
$$\forall N \in \nn, \forall f \in L^2(\rr^n), \quad \|\pi_Nf\|_{L^2(\rr^n)} \leq  Ce^{\frac{1}{2}N \ln(N+1)+CN}\|\pi_Nf\|_{L^2(\omega)}.$$
$(ii)$ If the measurable subset $\omega \subset \rr^n$ satisfies the condition \emph{(\ref{liminf})}, then there exists a positive constant $C=C(\omega)>1$ such that
$$\forall N \in \nn, \forall f \in L^2(\rr^n), \quad \|\pi_Nf\|_{L^2(\rr^n)} \leq  Ce^{CN}\|\pi_Nf\|_{L^2(\omega)}.$$
$(iii)$ If the measurable subset $\omega \subset \rr^n$ is $\gamma$-thick at scale $L>0$ in the sense defined in \emph{(\ref{thick1v})}, then there exist a positive constant $C=C(L,\gamma,n)>0$ depending on the dimension $n \geq 1$ and the parameters $\gamma, L>0$, and a universal positive constant $\kappa=\kappa(n)>0$ only depending on the dimension such that 
$$\forall N \in \nn, \forall f \in L^2(\rr^n), \quad \|\pi_Nf\|_{L^2(\rr^n)} \leq  C\Big(\frac{\kappa}{\gamma}\Big)^{\kappa L\sqrt{N}}\|\pi_Nf\|_{L^2(\omega)}.$$
\end{theorem}

\medskip

According to the above result, the control on the growth of the positive constant $C_N(\omega)>0$ with respect to the energy level for an arbitrary non-empty open subset $\omega$ of $\rr^n$, or when the measurable subset $\omega \subset \rr^n$ verifies condition (\ref{liminf}), is not sufficient to satisfy the estimates (\ref{eq7}) needed to obtain some results of null-controllability and observability for the parabolic equations associated to the class of hypoelliptic quadratic operators studied in Section~\ref{quadratic1}. As the one-dimensional harmonic heat equation is known from~\cite[Proposition~5.1]{miller1}, see also \cite{miller2}, to not be null-controllable, nor observable, in any time $T>0$ from a half-line and as the harmonic oscillator obviously belongs to the class of hypoelliptic quadratic operators studied in Section~\ref{quadratic1},     
we observe that spectral estimates of the type 
$$\exists 0<a<1, \exists C>1,  \forall N \in \nn, \forall f \in L^2(\rr^n), \quad \|\pi_Nf\|_{L^2(\rr^n)} \leq  Ce^{CN^a}\|\pi_Nf\|_{L^2(\omega)},$$
cannot hold for an arbitrary non-empty open subset $\omega$ of $\rr^n$, nor when the measurable subset $\omega \subset \rr^n$ verifies condition (\ref{liminf}),
since Theorem~\ref{Meta_thm_AdaptedLRmethod} together with (\ref{eq6}) would then imply the null-controlllability and the observability of the one-dimensional harmonic heat equation from a half-line. This would be in contradiction with the results of~\cite{miller1,miller2}.

On the other hand, when the measurable subset $\omega \subset \rr^n$ is $\gamma$-thick at scale $L>0$, the above spectral estimate $(iii)$ is an analogue for finite combinations of Hermite functions of the sharpened version of the Logvinenko-Sereda theorem proved by Kovrijkine in~\cite{thesis,Kovrijkine} with a similar dependence of the constant with respect to the parameters $0<\gamma \leq 1$ and $L>0$ as in (\ref{kov}). Notice that the growth in $\sqrt{N}$ is of the order of the square root of the largest eigenvalue of the harmonic oscillator $\mathcal{H}=-\Delta_x+|x|^2$ on the spectral vector subspace $\mathcal E_{N}$, whereas the growth in $R$ in (\ref{kov}) is also of order of the square root of the largest spectral value of the Laplace operator $-\Delta_x$ on the spectral vector subspace  
$$E_R=\{f \in L^2(\rr^n) : \textrm{ supp } \widehat{f} \subset [-R,R]^n\}.$$
This is in agreement with what is usually expected to hold for that type of spectral estimates, see~\cite{rousseau}.

The spectral estimate~$(i)$ for arbitrary non-empty open subsets is proved in Section~\ref{nonempty}. Its proof uses some estimates on Hermite functions together with the Remez inequality. The spectral estimate~$(ii)$ for measurable subsets satisfying condition (\ref{liminf}) is proved in Section~\ref{liminf1.3} and follows from similar arguments as the ones used in Section~\ref{nonempty}.
The spectral estimate~$(iii)$ for thick sets is proved in Section~\ref{thickproof}. 
\section{Proof of Theorem~\ref{th1}}

\subsection{Preliminary results}

\subsubsection{Hermite functions}\label{hermite_functions}
This section is devoted to set some notations and recall basic facts about Hermite functions.
The standard Hermite functions $(\phi_{k})_{k\geq 0}$ are defined for $x \in \rr$,
 \begin{equation}\label{defi}
 \phi_{k}(x)=\frac{(-1)^k}{\sqrt{2^k k!\sqrt{\pi}}} e^{\frac{x^2}{2}}\frac{d^k}{dx^k}(e^{-x^2})
 =\frac{1}{\sqrt{2^k k!\sqrt{\pi}}} \Bigl(x-\frac{d}{dx}\Bigr)^k(e^{-\frac{x^2}{2}})=\frac{ a_{+}^k \phi_{0}}{\sqrt{k!}},
\end{equation}
where $a_{+}$ is the creation operator
$$a_{+}=\frac{1}{\sqrt{2}}\Big(x-\frac{d}{dx}\Big).$$
The Hermite functions satisfy the identity 
\begin{equation}\label{sd1}
\forall k \geq 0,  \quad \widehat{\phi_{k}}=(-i)^k\sqrt{2\pi}\phi_k,
\end{equation}
when using the normalization of the Fourier transform (\ref{fourier}).
The $L^2$-adjoint of the creation operator is the annihilation operator
$$a_-=a_+^*=\frac{1}{\sqrt{2}}\Big(x+\frac{d}{dx}\Big).$$
The following identities hold
\begin{equation}\label{eq2ui}
[a_-,a_+]=1, \quad -\frac{d^2}{dx^2}+x^2=2a_+a_-+1,
\end{equation}
\begin{equation}\label{eq2}
\forall k \in \nn, \quad a_+ \phi_{k}=\sqrt{k+1} \phi_{k+1}, \qquad  \forall k \in \nn, \quad a_-\phi_{k}=\sqrt{k} \phi_{k-1} \ (=0\textrm{ if }k=0),
\end{equation}
\begin{equation}\label{eq2ui1}
\forall k \in \nn, \quad \Big(-\frac{d^2}{dx^2}+x^2\Big)\phi_{k}=(2k+1)\phi_{k},
\end{equation}
where $\nn$ denotes the set of non-negative integers.
The family $(\phi_{k})_{k\in \nn}$ is an orthonormal basis of $L^2(\R)$.
We set for $\alpha=(\alpha_{j})_{1\le j\le n}\in\N^n$, $x=(x_{j})_{1\le j\le n}\in \R^n,$
\begin{equation}\label{jk1}
\Phi_{\alpha}(x)=\prod_{j=1}^n\phi_{\alpha_j}(x_j).
\end{equation}
The family $(\Phi_{\alpha})_{\alpha \in \nn^n}$ is an orthonormal basis of $L^2(\R^n)$
composed of the eigenfunctions of the $n$-dimensional harmonic oscillator
\begin{equation}\label{6.harmo}
\mathcal{H}=-\Delta_x+|x|^2=\sum_{k\ge 0}(2k+n)\mathbb P_{k},\quad \text{Id}=\sum_{k \ge 0}\mathbb P_{k},
\end{equation}
where $\mathbb P_{k}$ is the orthogonal projection onto $\text{Span}_{\cc}
\{\Phi_{\alpha}\}_{\alpha\in \N^n,\val \alpha =k}$, with $\val \alpha=\alpha_{1}+\dots+\alpha_{n}$. 

The following estimates on Hermite functions are a key ingredient for the proof of the spectral inequalities $(i)$ and $(ii)$ in Theorem~\ref{th1}. This result was established by Bonami, Karoui and the second author in the proof of~\cite[Theorem~3.2]{BJK}, and is recalled here for the sake of completeness of the present work: 

\medskip

\begin{lemma}\label{lem:tail1}
The one-dimensional Hermite functions $(\phi_k)_{k \in \nn}$ defined in \emph{(\ref{defi})} satisfy the following estimates:
$$\forall k \in \nn, \forall a \geq \sqrt{2k+1}, \quad \int_{|x| \geq a}|\phi_k(x)|^2dx \leq \frac{2^{k+1}}{k!\sqrt{\pi}}a^{2k-1}e^{-a^2}.$$ 
\end{lemma}

\medskip

\begin{proof}
For any $k \in \nn$, the $k^{\textrm{th}}$ Hermite polynomial function
\begin{equation}\label{sz-1}
H_k(x)=(-1)^ke^{x^2}\Big(\frac{d}{dx}\Big)^{k}(e^{-x^2}),
\end{equation}
has degree $k$ and is an even (respectively odd) function when $k$ is an even (respectively odd) non-negative integer. The first Hermite polynomial functions are given by
\begin{equation}\label{sz}
H_0(x)=1, \qquad H_1(x)=2x, \qquad H_2(x)=4x^2-2.
\end{equation}
The $k^{\textrm{th}}$ Hermite polynomial function $H_k$ admits $k$ distinct real simple roots. More specifically, we recall from~\cite[Section~6.31]{szego} that the $k$ roots of $H_k$ denoted $-x_{[\frac{k}{2}],k}$, ..., $-x_{1,k}$, $x_{1,k}$, ..., $x_{[\frac{k}{2}],k}$, satisfy 
\begin{equation}\label{sz0}
-\sqrt{2k+1} \leq -x_{[\frac{k}{2}],k}<  ... <-x_{1,k}<0< x_{1,k} < ...<x_{[\frac{k}{2}],k} \leq \sqrt{2k+1},
\end{equation}
with $[\frac{k}{2}]$ the integer part of $\frac{k}{2}$, when $k \geq 2$ is an even positive integer. On the other hand, the $k$ roots of $H_k$ denoted $-x_{[\frac{k}{2}],k}$, ..., $-x_{1,k}$, $x_{0,k}$, $x_{1,k}$, ..., $x_{[\frac{k}{2}],k}$, satisfy 
\begin{equation}\label{sz1}
-\sqrt{2k+1} \leq -x_{[\frac{k}{2}],k}<  ... <-x_{1,k}<x_{0,k}=0< x_{1,k} < ...<x_{[\frac{k}{2}],k} \leq \sqrt{2k+1},
\end{equation}
when $k$ is an odd positive integer. We denote by $z_k$ the largest non-negative root of the $k^{\textrm{th}}$ Hermite polynomial function $H_k$, that is, with the above notations $z_k=x_{[\frac{k}{2}],k}$, when $k \geq 1$. Relabeling temporarily $(a_j)_{1 \leq j \leq k}$ the $k$ roots of $H_k$ such that 
$$a_1<a_2<...<a_k.$$
The classical formula 
\begin{equation}\label{sz6}
\forall k \in \nn^*, \quad H_k'(x)=2kH_{k-1}(x),
\end{equation}
see e.g.~\cite[Section~5.5]{szego}, together with Rolle's Theorem imply that $H_{k-1}$ admits exactly one root in each of the $k-1$ intervals $(a_j,a_{j+1})$, with $1 \leq j \leq k-1$, when $k \geq 2$. According to (\ref{sz}), (\ref{sz0}) and (\ref{sz1}), it implies in particular that for all $k \geq 1$,
\begin{equation}\label{sz2}
0=z_1<z_2<...<z_k \leq \sqrt{2k+1}.
\end{equation}
Next, we claim that 
\begin{equation}\label{sz3}
\forall k \geq 1, \forall |x| \geq z_k, \quad |H_k(x)| \leq 2^k|x|^k.
\end{equation}
To that end, we first observe that 
\begin{equation}\label{sz5}
\forall k \geq 1, \forall x \geq z_k, \quad H_k(x) \geq 0,
\end{equation}
since the leading coefficient of $H_k \in \rr[X]$ is given by $2^k>0$.
As the polynomial $H_k$ is an even or odd function, we notice from (\ref{sz5}) that it is actually sufficient to establish that  
\begin{equation}\label{sz4}
\forall k \geq 1, \forall x \geq z_k, \quad H_k(x) \leq 2^kx^k,
\end{equation}
to prove the claim. The estimates (\ref{sz4}) are proved by recurrence on $k \geq 1$. Indeed, we observe from (\ref{sz}) that 
$$\forall x \geq z_1=0, \quad H_1(x)=2x.$$
Let $k \geq 2$ such that the estimate (\ref{sz4}) is satisfied at rank $k-1$. It follows from (\ref{sz6}) for all $x \geq z_k$,
\begin{multline}\label{sz7}
H_k(x)=H_k(x)-H_k(z_k)=\int_{z_k}^xH_k'(t)dt=2k\int_{z_k}^xH_{k-1}(t)dt \\ 
\leq 2k\int_{z_k}^x2^{k-1}t^{k-1}dt=2^k(x^k-z_k^k) \leq 2^kx^k,
\end{multline}
since $0 \leq z_{k-1}<z_k$. This ends the proof of the claim (\ref{sz3}). We deduce from (\ref{sz}), (\ref{sz2}) and (\ref{sz3}) that 
\begin{equation}\label{sz8}
\forall k \in \nn, \forall |x| \geq \sqrt{2k+1}, \quad |H_k(x)| \leq 2^k|x|^k.
\end{equation}
It follows from (\ref{defi}), (\ref{sz-1}) and (\ref{sz8}) that 
\begin{equation}\label{sz9}
\forall k \in \nn, \forall |x| \geq \sqrt{2k+1}, \quad |\phi_{k}(x)|\leq\frac{2^{\frac{k}{2}}}{\sqrt{k!}\pi^{\frac{1}{4}}}|x|^k e^{-\frac{x^2}{2}}.
\end{equation}
We observe that 
\begin{equation}\label{sz10}
\forall a>0, \quad   \int_a^{+\infty}e^{-t^2}dt \leq a^{-1}e^{-\frac{a^2}{2}}\int_a^{+\infty}te^{-\frac{t^2}{2}}dt=a^{-1}e^{-a^2}
\end{equation}
and
\begin{equation}\label{sz11}
\forall \alpha>1, \forall a>\sqrt{\alpha-1}, \quad   \int_a^{+\infty}t^{\alpha}e^{-t^2}dt \leq a^{\alpha-1}e^{-\frac{a^2}{2}}\int_a^{+\infty}te^{-\frac{t^2}{2}}dt=a^{\alpha-1}e^{-a^2},
\end{equation}
as the function $(a,+\infty) \ni t \mapsto t^{\alpha-1}e^{-\frac{t^2}{2}} \in (0,+\infty)$ is decreasing on $(a,+\infty)$.
We deduce from (\ref{sz9}), (\ref{sz10}) and (\ref{sz11}) that 
\begin{multline}\label{sz12}
\forall k \in \nn, \forall a \geq \sqrt{2k+1}, \quad \int_{|x| \geq a}|\phi_{k}(x)|^2dx \leq \frac{2^{k}}{k!\pi^{\frac{1}{2}}}\int_{|x| \geq a}x^{2k} e^{-x^2}dx\\
=\frac{2^{k+1}}{k!\pi^{\frac{1}{2}}}\int_{x \geq a}x^{2k} e^{-x^2}dx \leq \frac{2^{k+1}}{k!\pi^{\frac{1}{2}}}a^{2k-1}e^{-a^2}.
\end{multline}
This ends the proof of Lemma~\ref{lem:tail1}.
\end{proof}

The following lemma is also instrumental in the proof of Theorem~\ref{th1}:  

\medskip

\begin{lemma}\label{lem:tail}
With $\mathcal E_{N}=\emph{\text{Span}}_{\cc}\{\Phi_{\alpha}\}_{\alpha\in \N^n, \val \alpha \leq N}$,
there exists a positive constant $c_n>0$ depending only on the dimension $n \geq 1$ such that
$$\forall N \in \nn, \forall f \in \mathcal E_{N}, \quad \int_{|x|\geq c_n\sqrt{N+1}}|f(x)|^2dx \leq \frac{1}{4}\norm{f}_{L^2(\R^n)}^2,$$
where $|\cdot|$ denotes the Euclidean norm on $\rr^n$.
\end{lemma}

\medskip

\begin{proof} Let $N \in \nn$. We deduce from Lemma~\ref{lem:tail1} and the Cauchy-Schwarz inequality that the one-dimensional Hermite functions $(\phi_{k})_{k\in \nn}$ satisfy for all $0 \leq k,l \leq N$, $a \geq \sqrt{2N+1}$,  
\begin{multline}\label{bona1}
\int_{|t|\geq a}|\phi_k(t)\phi_l(t)|dt
\leq \Big(\int_{|t| \geq a}|\phi_k(t)|^2dt\Big)^{\frac{1}{2}}\Big(\int_{|t| \geq a}|\phi_l(t)|^2dt\Big)^{\frac{1}{2}} \\ \leq 
\frac{2^{\frac{k+l}{2}+1}}{\sqrt{\pi} \sqrt{k!} \sqrt{l!}}a^{k+l-1}e^{-a^2}.
\end{multline}
In order to extend these estimates in the multi-dimensional setting, we first notice that for all $a>0$, $\alpha$, $\beta \in \nn^n$,  $|\alpha|, |\beta| \leq N$,  
\begin{equation}\label{lev1}
\int_{|x|\geq a}|\Phi_{\alpha}(x)\Phi_{\beta}(x)|dx \leq \sum_{j=1}^n\int_{|x_j|\geq \frac{a}{\sqrt{n}}}|\Phi_{\alpha}(x)\Phi_{\beta}(x)|dx.
\end{equation}
On the other hand, we notice from (\ref{bona1}) and (\ref{lev1}) that
\begin{align*}
& \ \int_{|x_j|\geq \frac{a}{\sqrt{n}}}|\Phi_{\alpha}(x)\Phi_{\beta}(x)|dx\\
=& \ \Big(\int_{|x_j|\geq \frac{a}{\sqrt{n}}}|\phi_{\alpha_j}(x_j)\phi_{\beta_j}(x_j)|dx_j\Big) \prod_{\substack{1 \leq k\leq n\\ k \neq j}}\Big(\int_{\rr}|\phi_{\alpha_k}(x_k)\phi_{\beta_k}(x_k)|dx_k\Big)
\\
\leq & \ \Big(\int_{|x_j|\geq \frac{a}{\sqrt{n}}}|\phi_{\alpha_j}(x_j)\phi_{\beta_j}(x_j)|dx_j\Big) \prod_{\substack{1 \leq k\leq n\\ k \neq j}}\|\phi_{\alpha_k}\|_{L^2(\rr)}\|\phi_{\beta_k}\|_{L^2(\rr)},
\end{align*}
implies that 
for all $a \geq \sqrt{n}\sqrt{2N+1}$, $\alpha$, $\beta \in \nn^n$,  $|\alpha|, |\beta| \leq N$,  
\begin{multline}\label{lev2}
\int_{|x|\geq a}|\Phi_{\alpha}(x)\Phi_{\beta}(x)|dx \leq \sum_{j=1}^n\int_{|x_j|\geq \frac{a}{\sqrt{n}}}|\phi_{\alpha_j}(x_j)\phi_{\beta_j}(x_j)|dx_j \\
\leq 2\sqrt{\frac{n}{\pi}}\frac{e^{-\frac{a^2}{n}}}{a}\sum_{j=1}^n\frac{1}{ \sqrt{\alpha_j!} \sqrt{\beta_j!}}\Big(\sqrt{\frac{2}{n}}a\Big)^{\alpha_j+\beta_j},
\end{multline}
since $(\phi_k)_{k \in \nn}$ is an orthonormal basis of $L^2(\rr)$.
For any $f=\sum_{|\alpha| \leq N} \gamma_{\alpha}\Phi_{\alpha} \in\mathcal{E}_N$ and $a\geq \sqrt{n}\sqrt{2N+1}$, we deduce from (\ref{lev2}) that 
\begin{multline}\label{lev3}
\int_{|x| \geq a}|f(x)|^2dx=\sum_{\substack{|\alpha| \leq N\\ |\beta| \leq N}} \gamma_{\alpha} \overline{\gamma_{\beta}} \int_{|x|\geq a}\Phi_{\alpha}(x)\overline{\Phi_{\beta}(x)}dx \\ \leq
\sum_{\substack{|\alpha| \leq N\\ |\beta| \leq N}} |\gamma_{\alpha}| |\gamma_{\beta}| \int_{|x|\geq a}|\Phi_{\alpha}(x)\Phi_{\beta}(x)|dx
\leq 2\sqrt{\frac{n}{\pi}}\frac{e^{-\frac{a^2}{n}}}{a} \sum_{\substack{|\alpha| \leq N,\ |\beta| \leq N\\
1 \leq j \leq n}}\frac{|\gamma_{\alpha}||\gamma_{\beta}|}{ \sqrt{\alpha_j!} \sqrt{\beta_j!}}\Big(\sqrt{\frac{2}{n}}a\Big)^{\alpha_j+\beta_j}.
\end{multline}
For any $\alpha=(\alpha_1,...,\alpha_n) \in \nn^n$, we denote $\alpha'=(\alpha_2,...,\alpha_n) \in \nn^{n-1}$ when $n \geq 2$. We observe that 
\begin{equation}\label{lev4}
\sum_{\substack{|\alpha| \leq N\\ |\beta| \leq N}}\frac{|\gamma_{\alpha}||\gamma_{\beta}|}{ \sqrt{\alpha_1!} \sqrt{\beta_1!}}\Big(\sqrt{\frac{2}{n}}a\Big)^{\alpha_1+\beta_1}=
\sum_{\substack{|\alpha'| \leq N\\ |\beta'| \leq N}}\Big(\sum_{\substack{0 \leq \alpha_1 \leq N-|\alpha'|\\ 0 \leq \beta_1 \leq N-|\beta'|}} \frac{|\gamma_{\alpha_1,\alpha'}||\gamma_{\beta_1,\beta'}|}{ \sqrt{\alpha_1!} \sqrt{\beta_1!}}\Big(\sqrt{\frac{2}{n}}a\Big)^{\alpha_1+\beta_1}\Big)
\end{equation}
and 
\begin{multline}\label{lev5}
\sum_{\substack{0 \leq \alpha_1 \leq N-|\alpha'|\\ 0 \leq \beta_1 \leq N-|\beta'|}} \frac{|\gamma_{\alpha_1,\alpha'}||\gamma_{\beta_1,\beta'}|}{ \sqrt{\alpha_1!} \sqrt{\beta_1!}}\Big(\sqrt{\frac{2}{n}}a\Big)^{\alpha_1+\beta_1} \\
\leq \Big(\sum_{\substack{0 \leq \alpha_1 \leq N-|\alpha'|\\ 0 \leq \beta_1 \leq N-|\beta'|}} | \gamma_{\alpha_1,\alpha'}|^2|\gamma_{\beta_1,\beta'}|^2 \Big)^{\frac{1}{2}}\Big(\sum_{\substack{0 \leq \alpha_1 \leq N-|\alpha'|\\ 0 \leq \beta_1 \leq N-|\beta'|}} \frac{(\frac{2a^2}{n})^{\alpha_1+\beta_1}}{\alpha_1! \beta_1!}\Big)^{\frac{1}{2}},
\end{multline}
thanks to the Cauchy-Schwarz inequality. On the other hand, we notice that 
\begin{equation}\label{lev6}
\Big(\sum_{\substack{0 \leq \alpha_1 \leq N-|\alpha'|\\ 0 \leq \beta_1 \leq N-|\beta'|}} \frac{(\frac{2a^2}{n})^{\alpha_1+\beta_1}}{\alpha_1! \beta_1!}\Big)^{\frac{1}{2}} \leq 4^{N}\Big(\sum_{\substack{0 \leq \alpha_1 \leq N-|\alpha'|\\ 0 \leq \beta_1 \leq N-|\beta'|}} \frac{(\frac{a^2}{2n})^{\alpha_1+\beta_1}}{\alpha_1! \beta_1!}\Big)^{\frac{1}{2}}\leq 4^{N}e^{\frac{a^2}{2n}}.
\end{equation}
It follows from (\ref{lev4}), (\ref{lev5}) and (\ref{lev6}) that 
\begin{equation}\label{lev7}
\sum_{\substack{|\alpha| \leq N\\ |\beta| \leq N}}\frac{|\gamma_{\alpha}||\gamma_{\beta}|}{ \sqrt{\alpha_1!} \sqrt{\beta_1!}}\Big(\sqrt{\frac{2}{n}}a\Big)^{\alpha_1+\beta_1} \leq 4^{N}e^{\frac{a^2}{2n}}
\sum_{\substack{|\alpha'| \leq N\\ |\beta'| \leq N}}\Big(\sum_{\substack{0 \leq \alpha_1 \leq N-|\alpha'|\\ 0 \leq \beta_1 \leq N-|\beta'|}} | \gamma_{\alpha_1,\alpha'}|^2|\gamma_{\beta_1,\beta'}|^2 \Big)^{\frac{1}{2}}.
\end{equation}
The Cauchy-Schwarz inequality implies that 
\begin{multline}\label{lev8}
\sum_{\substack{|\alpha'| \leq N\\ |\beta'| \leq N}}\Big(\sum_{\substack{0 \leq \alpha_1 \leq N-|\alpha'|\\ 0 \leq \beta_1 \leq N-|\beta'|}} | \gamma_{\alpha_1,\alpha'}|^2|\gamma_{\beta_1,\beta'}|^2 \Big)^{\frac{1}{2}} \\ \leq \Big(\sum_{\substack{|\alpha'| \leq N\\ |\beta'| \leq N}}\Big(\sum_{\substack{0 \leq \alpha_1 \leq N-|\alpha'|\\ 0 \leq \beta_1 \leq N-|\beta'|}} | \gamma_{\alpha_1,\alpha'}|^2|\gamma_{\beta_1,\beta'}|^2 \Big)\Big)^{\frac{1}{2}}\Big(\sum_{\substack{|\alpha'| \leq N\\ |\beta'| \leq N}}1\Big)^{\frac{1}{2}}.
\end{multline}
By using that  the family $(\Phi_{\alpha})_{\alpha \in \nn^n}$ is an orthonormal basis of $L^2(\R^n)$ and that the number of solutions to the equation $\alpha_2+...+\alpha_{n}=k$, with $k \geq 0$, $n \geq 2$ and unknown $\alpha'=(\alpha_2,...,\alpha_n) \in \nn^{n-1}$, is given by $\binom{k+n-2}{n-2}$, we deduce from (\ref{lev8}) that 
\begin{multline}\label{lev9}
\sum_{\substack{|\alpha'| \leq N\\ |\beta'| \leq N}}\Big(\sum_{\substack{0 \leq \alpha_1 \leq N-|\alpha'|\\ 0 \leq \beta_1 \leq N-|\beta'|}} | \gamma_{\alpha_1,\alpha'}|^2|\gamma_{\beta_1,\beta'}|^2 \Big)^{\frac{1}{2}}  \leq \Big(\sum_{|\alpha| \leq N}| \gamma_{\alpha}|^2\Big)\Big(\sum_{|\alpha'| \leq N}1\Big)\\
=\Big(\sum_{k=0}^N\binom{k+n-2}{n-2}\Big)\|f\|_{L^2(\rr^n)}^2 \leq 2^{n-2}\Big(\sum_{k=0}^N2^{k}\Big)\|f\|_{L^2(\rr^n)}^2 \leq 2^{N+n-1}\|f\|_{L^2(\rr^n)}^2,
\end{multline}
since $\binom{k+n-2}{n-2} \leq \sum_{j=0}^{k+n-2}\binom{k+n-2}{j} = 2^{k+n-2}$.
It follows from (\ref{lev7}) and (\ref{lev9}) that 
\begin{equation}\label{lev10}
\sum_{\substack{|\alpha| \leq N\\ |\beta| \leq N}}\frac{|\gamma_{\alpha}||\gamma_{\beta}|}{ \sqrt{\alpha_1!} \sqrt{\beta_1!}}\Big(\sqrt{\frac{2}{n}}a\Big)^{\alpha_1+\beta_1} \leq 2^{n-1}8^{N}e^{\frac{a^2}{2n}}\|f\|_{L^2(\rr^n)}^2,
\end{equation}
when $n \geq 2$. Notice that the very same estimate holds true as well in the one-dimensional case $n=1$. 
We deduce from (\ref{lev3}) and (\ref{lev10}) that for all $N \in \nn$, $f \in \mathcal{E}_N$, $a\geq \sqrt{n}\sqrt{2N+1}$,     
\begin{equation}\label{lev11}
\int_{|x| \geq a}|f(x)|^2dx
\leq \frac{2^nn^{\frac{3}{2}}}{\sqrt{\pi}}\frac{e^{-\frac{a^2}{2n}}}{a}8^{N}\|f\|_{L^2(\rr^n)}^2 .
\end{equation}
It follows from (\ref{lev11}) that there exists a positive constant $c_n>0$ depending only on the dimension $n \geq 1$ such that
$$\forall N \in \nn, \forall f \in \mathcal E_{N}, \quad \int_{|x|\geq c_n\sqrt{N+1}}|f(x)|^2dx \leq \frac{1}{4}\norm{f}_{L^2(\R^n)}^2.$$ 
This ends the proof of Lemma~\ref{lem:tail}.

\end{proof}

\subsubsection{Bernstein type estimates for Hermite functions}\label{weighted}
This section is devoted to the proof of the following Bernstein type estimates for Hermite functions:

\medskip

\begin{proposition}\label{prop1}
With $\mathcal E_{N}=\emph{\text{Span}}_{\cc}\{\Phi_{\alpha}\}_{\alpha\in \N^n, \val \alpha \leq N}$, finite combinations of Hermite functions satisfy the following estimates:
\begin{multline*}
(i) \quad \forall N \in \nn, \forall f \in \mathcal E_{N}, \forall 0<\delta \leq 1, \forall \beta \in \nn^n, \\ \|\partial_x^{\beta}f\|_{L^2(\rr^n)}\leq e^{\frac{e}{2\delta^2}}(2\delta)^{|\beta|}|\beta|!e^{\delta^{-1}\sqrt{N}}\|f\|_{L^2(\rr^n)}.
\end{multline*}
\begin{multline*}
(ii) \quad \forall N \in \nn, \forall f \in \mathcal E_{N}, \forall 0<\delta<\frac{1}{32n}, \forall \beta \in \nn^n, \\ \|e^{\delta |x|^2}\partial_x^{\beta}f\|_{L^2(\rr^n)} +\|e^{\delta |D_x|^2}x^{\beta}f\|_{L^2(\rr^n)}
\leq \frac{2^{n}}{1-32n \delta}2^{\frac{N}{2}}2^{\frac{3}{2}|\beta|}\sqrt{|\beta|!}\|f\|_{L^2(\rr^n)}.
\end{multline*}
\end{proposition}

\medskip

\begin{proof}
We notice that 
\begin{equation}\label{eq1}
x_j=\frac{1}{\sqrt{2}}(a_{j,+}+a_{j,-}), \quad \partial_{x_j}=\frac{1}{\sqrt{2}}(a_{j,-}-a_{j,+}),
\end{equation}
with 
$$a_{j,+}=\frac{1}{\sqrt{2}}(x_j-\partial_{x_j}), \quad a_{j,-}=\frac{1}{\sqrt{2}}(x_j+\partial_{x_j}).$$
By denoting $(e_j)_{1 \leq j \leq n}$ the canonical basis of $\rr^n$, we obtain from (\ref{eq2}) and (\ref{eq1}) that for all $N \in \nn$, $f \in \mathcal E_{N}$,
\begin{align*}
& \ \|a_{j,+}f\|_{L^2(\rr^n)}^2=\Big\|a_{j,+}
\Big(\sum_{|\alpha| \leq N}\langle f,\Phi_{\alpha}\rangle_{L^2}\Phi_{\alpha}\Big)\Big\|_{L^2(\rr^n)}^2\\
=& \ \Big\|\sum_{|\alpha| \leq N}\sqrt{\alpha_j+1}\langle f,\Phi_{\alpha}\rangle_{L^2}\Phi_{\alpha+e_j}\Big\|_{L^2(\rr^n)}^2
=\sum_{|\alpha| \leq N}(\alpha_j+1)|\langle f,\Phi_{\alpha}\rangle_{L^2}|^2\\
\leq & \ (N+1)\sum_{|\alpha| \leq N}|\langle f,\Phi_{\alpha}\rangle_{L^2}|^2=(N+1)\|f\|_{L^2(\rr^n)}^2
\end{align*}
and
\begin{align*}
& \ \|a_{j,-}f\|_{L^2(\rr^n)}^2=\Big\|a_{j,-}
\Big(\sum_{|\alpha| \leq N}\langle f,\Phi_{\alpha}\rangle_{L^2}\Phi_{\alpha}\Big)\Big\|_{L^2(\rr^n)}^2\\
=& \ \Big\|\sum_{|\alpha| \leq N}\sqrt{\alpha_j}\langle f,\Phi_{\alpha}\rangle_{L^2}\Phi_{\alpha-e_j}\Big\|_{L^2(\rr^n)}^2
=\sum_{|\alpha| \leq N}\alpha_j|\langle f,\Phi_{\alpha}\rangle_{L^2}|^2\\
\leq & \ N\sum_{|\alpha| \leq N}|\langle f,\Phi_{\alpha}\rangle_{L^2}|^2=N\|f\|_{L^2(\rr^n)}^2.
\end{align*}
It follows that for all $N \in \nn$, $f \in \mathcal E_{N}$,
\begin{equation}\label{a1}
\|x_jf\|_{L^2(\rr^n)}\leq \frac{1}{\sqrt{2}}( \|a_{j,+}f\|_{L^2(\rr^n)}+ \|a_{j,-}f\|_{L^2(\rr^n)}) \leq \sqrt{2N+2}\|f\|_{L^2(\rr^n)}
\end{equation}
and
\begin{equation}\label{a2}
\|\partial_{x_j}f\|_{L^2(\rr^n)}\leq \frac{1}{\sqrt{2}}( \|a_{j,+}f\|_{L^2(\rr^n)}+ \|a_{j,-}f\|_{L^2(\rr^n)}) \leq \sqrt{2N+2}\|f\|_{L^2(\rr^n)}.
\end{equation}
We notice from (\ref{eq2}) and (\ref{eq1}) that
$$\forall N \in \nn, \forall f \in \mathcal E_{N}, \forall \alpha, \beta \in \nn^n, \quad x^{\alpha}\partial_x^{\beta}f \in  \mathcal E_{N+|\alpha|+|\beta|},$$
with $x^{\alpha}=x_1^{\alpha_1}...x_n^{\alpha_n}$ and $\partial_x^{\beta}=\partial_{x_1}^{\beta_1}...\partial_{x_n}^{\beta_n}$.
We deduce from (\ref{a1}) that for all $N \in \nn$, $f \in \mathcal E_{N}$, and $\alpha, \beta \in \nn^n$, with $\alpha_1 \geq 1$,
$$\|x^{\alpha}\partial_x^{\beta}f\|_{L^2(\rr^n)}=\|x_1(\underbrace{x^{\alpha-e_1}\partial_x^{\beta}f}_{\in \mathcal E_{N+|\alpha|+|\beta|-1}})\|_{L^2(\rr^n)}\leq
\sqrt{2}\sqrt{N+|\alpha|+|\beta|}\|x^{\alpha-e_1}\partial_x^{\beta}f\|_{L^2(\rr^n)}.$$
By iterating the previous estimates, we readily obtain from (\ref{a1}) and (\ref{a2}) that for all $N \in \nn$, $f \in \mathcal E_{N}$, $\alpha, \beta \in \nn^n$, 
\begin{equation}\label{gh0}
\|x^{\alpha}\partial_x^{\beta}f\|_{L^2(\rr^n)}\leq 2^{\frac{|\alpha|+|\beta|}{2}}\sqrt{\frac{(N+|\alpha|+|\beta|)!}{N!}}\|f\|_{L^2(\rr^n)}.
\end{equation}
We recall the following basic estimates, 
\begin{equation}\label{rod1}
\forall k \in \nn^*, \ k^k \leq e^k k!, \qquad \forall t,A>0, \ t^A \leq A^Ae^{t-A}, \qquad \forall t>0, \forall k \in \nn, \ t^k \leq e^t k!,
\end{equation}
see e.g.~\cite{rodino} (formulas (0.3.12) and (0.3.14)). 
Let $0<\delta \leq 1$ be a positive constant. When $N \leq |\alpha|+|\beta|$, we deduce from (\ref{rod1}) that 
\begin{multline}\label{gh1}
2^{\frac{|\alpha|+|\beta|}{2}}\sqrt{\frac{(N+|\alpha|+|\beta|)!}{N!}} \leq 2^{\frac{|\alpha|+|\beta|}{2}}(N+|\alpha|+|\beta|)^{\frac{|\alpha|+|\beta|}{2}}\leq 2^{|\alpha|+|\beta|}(|\alpha|+|\beta|)^{\frac{|\alpha|+|\beta|}{2}}\\
\leq (2\sqrt{e})^{|\alpha|+|\beta|}\sqrt{(|\alpha|+|\beta|)!}=(2\sqrt{e})^{|\alpha|+|\beta|}\frac{(|\alpha|+|\beta|)!}{\sqrt{(|\alpha|+|\beta|)!}}\leq
e^{\frac{e}{2\delta^2}}(2\delta)^{|\alpha|+|\beta|}(|\alpha|+|\beta|)!.
\end{multline}
On the other hand, when $N \geq |\alpha|+|\beta|$, we deduce from (\ref{rod1}) that 
\begin{align}\label{gh2}
& \ 2^{\frac{|\alpha|+|\beta|}{2}}\sqrt{\frac{(N+|\alpha|+|\beta|)!}{N!}} \leq 2^{\frac{|\alpha|+|\beta|}{2}}(N+|\alpha|+|\beta|)^{\frac{|\alpha|+|\beta|}{2}}\\ \notag
\leq & \ (2\delta)^{|\alpha|+|\beta|}(\delta^{-1} \sqrt{N})^{|\alpha|+|\beta|}
\leq (2\delta)^{|\alpha|+|\beta|}(|\alpha|+|\beta|)^{|\alpha|+|\beta|}e^{\delta^{-1}\sqrt{N}-|\alpha|-|\beta|}\\ \notag
 \leq & \  (2\delta)^{|\alpha|+|\beta|}(|\alpha|+|\beta|)!e^{\delta^{-1}\sqrt{N}}.
\end{align}
It follows from (\ref{gh0}), (\ref{gh1}) and (\ref{gh2}) that for all $N \in \nn$, $f \in \mathcal E_{N}$, $\alpha, \beta \in \nn^n$, 
\begin{equation}\label{gh3}
\|x^{\alpha}\partial_x^{\beta}f\|_{L^2(\rr^n)}\leq e^{\frac{e}{2\delta^2}} (2\delta)^{|\alpha|+|\beta|}(|\alpha|+|\beta|)!e^{\delta^{-1}\sqrt{N}}\|f\|_{L^2(\rr^n)}.
\end{equation}
This provides in particular the following Bernstein type estimates
\begin{multline}\label{bernstein1}
\forall N \in \nn, \forall f \in \mathcal E_{N}, \forall 0<\delta \leq 1, \forall \beta \in \nn^n, \\  \|\partial_x^{\beta}f\|_{L^2(\rr^n)}\leq e^{\frac{e}{2\delta^2}}(2\delta)^{|\beta|}|\beta|!e^{\delta^{-1}\sqrt{N}}\|f\|_{L^2(\rr^n)}.
\end{multline}
On the other hand, we deduce from (\ref{gh0}) that for all $N \in \nn$, $f \in \mathcal E_{N}$, $\alpha, \beta \in \nn^n$, 
\begin{multline}\label{jk2}
\|x^{\alpha}\partial_x^{\beta}f\|_{L^2(\rr^n)}\leq 2^{\frac{|\alpha|+|\beta|}{2}}\sqrt{\frac{(N+|\alpha|+|\beta|)!}{N!}}\|f\|_{L^2(\rr^n)}\\
\leq 2^{\frac{N}{2}}2^{|\alpha|+|\beta|}\sqrt{(|\alpha|+|\beta|)!}\|f\|_{L^2(\rr^n)},
\end{multline}
since
$$\frac{(k_1+k_2)!}{k_1!k_2!}=\binom{k_1+k_2}{k_1} \leq \sum_{j=0}^{k_1+k_2}\binom{k_1+k_2}{j}=2^{k_1+k_2}.$$
We observe from (\ref{jk2}) that for all $N \in \nn$, $f \in \mathcal E_{N}$, $\delta>0$, $\alpha, \beta \in \nn^n$,
\begin{multline}\label{esti1}
\Big\|\frac{\delta^{|\alpha|}x^{2\alpha}}{\alpha!}\partial_x^{\beta}f\Big\|_{L^2(\rr^n)}\leq \frac{2^{\frac{N}{2}}\delta^{|\alpha|}2^{2|\alpha|+|\beta|}}{\alpha!}\sqrt{(2|\alpha|+|\beta|)!}\|f\|_{L^2(\rr^n)}\\
\leq 2^{\frac{N}{2}}\delta^{|\alpha|}2^{4|\alpha|+\frac{3}{2}|\beta|}\frac{|\alpha|!}{\alpha!}\sqrt{|\beta|!}\|f\|_{L^2(\rr^n)}
\leq 2^{\frac{N}{2}}(16n \delta)^{|\alpha|}2^{\frac{3}{2}|\beta|}\sqrt{|\beta|!}\|f\|_{L^2(\rr^n)},
\end{multline}
since
$$(2|\alpha|+|\beta|)! \leq 2^{2|\alpha|+|\beta|}(2|\alpha|)!|\beta|! \leq 2^{4|\alpha|+|\beta|}(|\alpha|!)^2|\beta|!$$
and
\begin{equation}\label{esti3}
|\alpha|! \leq n^{|\alpha|}\alpha!.
\end{equation}
The last estimate is a direct consequence of the generalized Newton formula
$$\forall x=(x_1,...,x_n) \in \rr^n, \forall N \in \nn, \quad \Big(\sum_{j=1}^nx_j\Big)^N=\sum_{\alpha \in \nn^n, |\alpha|=N}\frac{N!}{\alpha!}x^{\alpha}.$$
By using that the number of solutions to the equation $\alpha_1+...+\alpha_n=k$, with $k \geq 0$, $n \geq 1$ and unknown $\alpha=(\alpha_1,...,\alpha_n) \in \nn^n$, is given by $\binom{k+n-1}{n-1}$,
it follows from (\ref{esti1}) that for all $N \in \nn$, $f \in \mathcal E_{N}$, $0<\delta<\frac{1}{32n}$, $\beta \in \nn^n$,
\begin{align}\label{jk4}
\|e^{\delta |x|^2}\partial_x^{\beta}f\|_{L^2(\rr^n)} \leq & \  \sum_{\alpha \in \nn^n}\Big\|\frac{\delta^{|\alpha|}x^{2\alpha}}{\alpha!}\partial_x^{\beta}f\Big\|_{L^2(\rr^n)}
\\ \notag
\leq & \ 2^{\frac{N}{2}}\Big(\sum_{\alpha \in \nn^n}(16n \delta)^{|\alpha|}\Big)2^{\frac{3}{2}|\beta|}\sqrt{|\beta|!}\|f\|_{L^2(\rr^n)}\\ \notag
= & \ 2^{\frac{N}{2}}\Big(\sum_{k=0}^{+\infty}\binom{k+n-1}{n-1}(16n \delta)^{k}\Big)2^{\frac{3}{2}|\beta|}\sqrt{|\beta|!}\|f\|_{L^2(\rr^n)} \\ \notag
\leq & \ \frac{2^{n-1}}{1-32n \delta}2^{\frac{N}{2}}2^{\frac{3}{2}|\beta|}\sqrt{|\beta|!}\|f\|_{L^2(\rr^n)},
\end{align}
since $\binom{k+n-1}{n-1} \leq \sum_{j=0}^{k+n-1}\binom{k+n-1}{j}= 2^{k+n-1}$. By noticing from (\ref{sd1}) that $f \in \mathcal E_{N}$ if and only if $\widehat{f} \in \mathcal E_{N}$, we deduce from the Parseval formula and (\ref{jk4}) that for all $N \in \nn$, $f \in \mathcal E_{N}$, $0<\delta<\frac{1}{32n}$, $\beta \in \nn^n$,
\begin{multline}\label{jk5}
\|e^{\delta |D_x|^2}x^{\beta}f\|_{L^2(\rr^n)}=\frac{1}{(2\pi)^{\frac{n}{2}}}\|e^{\delta |\xi|^2}\partial_{\xi}^{\beta}\widehat{f}\|_{L^2(\rr^n)} \\ \leq 
 \frac{1}{(2\pi)^{\frac{n}{2}}}\frac{2^{n-1}}{1-32n \delta}2^{\frac{N}{2}}2^{\frac{3}{2}|\beta|}\sqrt{|\beta|!}\|\widehat{f}\|_{L^2(\rr^n)}=
\frac{2^{n-1}}{1-32n \delta}2^{\frac{N}{2}}2^{\frac{3}{2}|\beta|}\sqrt{|\beta|!}\|f\|_{L^2(\rr^n)}.
\end{multline}
This ends the proof of Proposition~\ref{prop1}.
\end{proof}

\subsection{Proofs of the uncertainty principles for Hermite functions}

This section is devoted to the proof of Theorem~\ref{th1}.

\subsubsection{Case when the control subset is a non-empty open set}\label{nonempty}
Let $\omega \subset \rr^n$ be a non-empty open set. There exist $x_0 \in \rr^n$ and $r>0$ such that the control subset $\omega$ contains the following open Euclidean ball 
\begin{equation}\label{sd3}
B(x_0,r) \subset \omega.
\end{equation}
We recall from (\ref{spec}) that 
\begin{equation}\label{spec1}
\forall N \in \nn, \exists C_N(\omega)>0, \forall f \in \mathcal E_{N}, \quad \|f\|_{L^2(\rr^n)} \leq C_N(\omega)\|f\|_{L^2(\omega)},
\end{equation} 
with $\mathcal E_{N}=\pi_N(L^2(\rr^n))$.
On the other hand, it follows from Lemma~\ref{lem:tail} that 
\begin{equation}\label{we1}
\forall N \in \nn, \forall f \in \mathcal E_{N}, \quad \|f\|_{L^2(\rr^n)} \leq \frac{2}{\sqrt{3}}\|f\|_{L^2(B(0,c_n\sqrt{N+1}))}.
\end{equation}
Let $N \in \nn$ and $f \in \mathcal E_{N}$. According to (\ref{defi}) and (\ref{jk1}), there exists a complex polynomial function $P \in \cc[X_1,...,X_n]$ of degree at most $N$ such that 
\begin{equation}\label{we2}
\forall x \in \rr^n, \quad f(x)=P(x)e^{-\frac{|x|^2}{2}}.
\end{equation}
We observe from (\ref{we1}) and (\ref{we2}) that 
\begin{equation}\label{we3}
\|f\|_{L^2(\rr^n)}^2 \leq \frac{4}{3}\int_{B(0,c_n\sqrt{N+1})}|P(x)|^2e^{-|x|^2}dx \leq \frac{4}{3}\|P\|_{L^2(B(0,c_n\sqrt{N+1}))}^2
\end{equation}
and
\begin{equation}\label{we4}
\|P\|_{L^2(B(x_0,r))}^2=\int_{B(x_0,r)}|P(x)|^2e^{-|x|^2}e^{|x|^2}dx \leq e^{(|x_0|+r)^2}\|f\|_{L^2(B(x_0,r))}^2.
\end{equation}
We aim at deriving an estimate of the term $\|P\|_{L^2(B(0,c_n\sqrt{N+1}))}$ by $\|P\|_{L^2(B(x_0,r))}$ when $N\gg 1$ is sufficiently large. Let $N$ be an integer such that $c_n\sqrt{N+1} >2|x_0|+r$. It implies the inclusion $B(x_0,r) \subset B(0,c_n\sqrt{N+1})$.
To that end, we may assume that $P$ is a non-zero polynomial function. 
By using polar coordinates centered at $x_0$, we notice that 
$$B(x_0,r)=\{x_0+t\sigma\,: 0\leq t < r,\ \sigma\in\mathbb{S}^{n-1}\}$$
and
\begin{equation}\label{we5}
\|P\|_{L^2(B(x_0,r))}^2=\int_{\mathbb{S}^{n-1}}\int_0^r|P(x_0+t\sigma)|^2t^{n-1}dtd\sigma.
\end{equation}
As $c_n\sqrt{N+1} >2|x_0|+r$, we notice that there exists a continuous function $\rho_N : \mathbb{S}^{n-1} \rightarrow (0,+\infty)$ such that 
\begin{equation}\label{we6}
B(0,c_n\sqrt{N+1})=\{x_0+t\sigma\ : 0\leq t< \rho_N(\sigma),\ \sigma\in\mathbb{S}^{n-1}\}
\end{equation}
and
\begin{equation}\label{we7}
\forall \sigma \in \mathbb{S}^{n-1}, \quad 0<|x_0|+r<c_n\sqrt{N+1}-|x_0|<\rho_N(\sigma)<c_n\sqrt{N+1}+|x_0|.
\end{equation}
It follows from (\ref{we6}) and (\ref{we7}) that 
\begin{multline}\label{we7.1}
\|P\|_{L^2(B(0,c_n\sqrt{N+1})\setminus B(x_0,\frac{r}{2}))}^2=\int_{\mathbb{S}^{n-1}}\int_{\frac{r}{2}}^{\rho_N(\sigma)}|P(x_0+t\sigma)|^2t^{n-1}dtd\sigma\\
\leq (c_n\sqrt{N+1}+|x_0|)^{n-1}\int_{\mathbb{S}^{n-1}}\int_{\frac{r}{2}}^{\rho_N(\sigma)}|P(x_0+t\sigma)|^2dtd\sigma.
\end{multline}
By noticing that 
$$t\to P\Big(x_0+(\frac{\rho_N(\sigma)}{2}+\frac{r}{4})\sigma+t\sigma\Big),$$
is a polynomial function of degree at most $N$, we deduce from (\ref{we7}) and Lemma~\ref{remez-1} used in the one-dimensional case $n=1$ that
\begin{align}\label{we8}
& \ \int_{\frac{r}{2}}^{\rho_N(\sigma)}|P(x_0+t\sigma)|^2dt =\int_{-(\frac{\rho_N(\sigma)}{2}-\frac{r}{4})}^{\frac{\rho_N(\sigma)}{2}-\frac{r}{4}}
\Big|P\Big(x_0+\Big(\frac{\rho_N(\sigma)}{2}+\frac{r}{4}\Big)\sigma+t\sigma\Big)\Big|^2dt\\ \notag
\leq & \ \frac{2^{4N+2}}{3}\frac{4(\rho_N(\sigma)-\frac{r}{2})}{\frac{r}{2}}
\left(\frac{2-\frac{\frac{r}{2}}{4(\rho_N(\sigma)-\frac{r}{2})}}{\frac{\frac{r}{2}}{4(\rho_N(\sigma)-\frac{r}{2})}}\right)^{2N}
\int_{-(\frac{\rho_N(\sigma)}{2}-\frac{r}{4})}^{\frac{-\rho_N(\sigma)}{2}+\frac{3r}{4}}
\Big|P\Big(x_0+\Big(\frac{\rho_N(\sigma)}{2}+\frac{r}{4}\Big)\sigma+t\sigma\Big)\Big|^2dt\\ \notag
\leq & \ \frac{2^{4N+2}}{3}\frac{4(\rho_N(\sigma)-\frac{r}{2})}{\frac{r}{2}}
\left(\frac{2-\frac{\frac{r}{2}}{4(\rho_N(\sigma)-\frac{r}{2})}}{\frac{\frac{r}{2}}{4(\rho_N(\sigma)-\frac{r}{2})}}\right)^{2N}
\int_{\frac{r}{2}}^{r}|P(x_0+t\sigma)|^2dt\\ \notag
\leq & \ \frac{2^{12N+n+4}}{3r^{2N+n}}\Big(c_n\sqrt{N+1}+|x_0|-\frac{r}{2}\Big)^{2N+1}\int_{\frac{r}{2}}^{r}|P(x_0+t\sigma)|^2t^{n-1}dt.
\end{align}
It follows from (\ref{we7.1}) and (\ref{we8}) that
\begin{multline}\label{we9}
\|P\|_{L^2(B(0,c_n\sqrt{N+1})\setminus B(x_0,\frac{r}{2}))}^2
\leq (c_n\sqrt{N+1}+|x_0|)^{n-1}\\ \times 
\frac{2^{12N+n+4}}{3r^{2N+n}}\Big(c_n\sqrt{N+1}+|x_0|-\frac{r}{2}\Big)^{2N+1}\int_{\mathbb{S}^{n-1}}\int_{\frac{r}{2}}^{r}|P(x_0+t\sigma)|^2t^{n-1}dt,
\end{multline}
implying that 
\begin{multline}\label{we10}
\|P\|_{L^2(B(0,c_n\sqrt{N+1}))}^2 
\leq \Big(1+(c_n\sqrt{N+1}+|x_0|)^{n-1}\\ \times \frac{2^{12N+n+4}}{3r^{2N+n}}\Big(c_n\sqrt{N+1}+|x_0|-\frac{r}{2}\Big)^{2N+1}\|P\|_{L^2(B(x_0,r))}^2,
\end{multline}
thanks to (\ref{we5}).
We deduce from (\ref{we10}) that there exists a positive constant $C=C(x_0,r,n)>1$ independent on the parameter $N$ such that 
\begin{equation}\label{we11}
\|P\|_{L^2(B(0,c_n\sqrt{N+1}))} 
\leq Ce^{\frac{1}{2}N\ln(N+1)+CN}\|P\|_{L^2(B(x_0,r))}.
\end{equation}
It follows from (\ref{we3}), (\ref{we4}) and (\ref{we11}) that for all $N \in \nn$ such that $c_n\sqrt{N+1} >2|x_0|+r$ and for all $f \in \mathcal{E}_N$,
\begin{equation}\label{we12}
\|f\|_{L^2(\rr^n)} \leq  \frac{2}{\sqrt{3}}Ce^{\frac{1}{2}(|x_0|+r)^2}e^{\frac{1}{2}N\ln(N+1)+CN}\|f\|_{L^2(B(x_0,r))}.
\end{equation}
The two estimates (\ref{spec1}) and (\ref{we12}) allow to prove the assertion $(i)$ in Theorem~\ref{th1}.

\subsubsection{Case when the control subset is a measurable set satisfying condition (\ref{liminf})}\label{liminf1.3}
Let $\omega \subset \rr^n$ be a measurable subset satisfying the condition 
\begin{equation}\label{liminf1}
\liminf_{R \to +\infty}\frac{|\omega \cap B(0,R)|}{|B(0,R)|}=\lim_{R \to +\infty}\Big(\inf_{r \geq R}\frac{|\omega \cap B(0,r)|}{|B(0,r)|}\Big)>0.
\end{equation}
It follows that there exist some positive constants $R_0>0$ and $\delta>0$ such that 
\begin{equation}\label{liminf2}
\forall R\geq R_0, \quad  \frac{|\omega \cap B(0,R)|}{|B(0,R)|} \geq \delta>0.
\end{equation}
We recall from (\ref{spec}) that 
\begin{equation}\label{spec1x}
\forall N \in \nn, \exists C_N(\omega)>0, \forall f \in \mathcal E_{N}, \quad \|f\|_{L^2(\rr^n)} \leq C_N(\omega)\|f\|_{L^2(\omega)}
\end{equation} 
and as in the above section, it follows from Lemma~\ref{lem:tail} that 
\begin{equation}\label{we1x}
\forall N \in \nn, \forall f \in \mathcal E_{N}, \quad \|f\|_{L^2(\rr^n)} \leq \frac{2}{\sqrt{3}}\|f\|_{L^2(B(0,c_n\sqrt{N+1}))},
\end{equation}
where $c_n>0$ is a positive constant depending only on the dimension $n \geq 1$.
Let $N \in \nn$ be an integer satisfying $c_n\sqrt{N+1} \geq R_0$ and $f \in \mathcal E_{N}$. It follows from (\ref{liminf2}) that 
\begin{equation}\label{dfx}
|\omega \cap B(0,c_n\sqrt{N+1})| \geq \delta |B(0,c_n\sqrt{N+1})|>0.
\end{equation}
According to (\ref{defi}) and (\ref{jk1}), there exists a complex polynomial function $P \in \cc[X_1,...,X_n]$ of degree at most $N$ such that 
\begin{equation}\label{we2x}
\forall x \in \rr^n, \quad f(x)=P(x)e^{-\frac{|x|^2}{2}}.
\end{equation}
We observe from (\ref{we1x}) and (\ref{we2x}) that 
\begin{equation}\label{we3x}
\|f\|_{L^2(\rr^n)}^2 \leq \frac{4}{3}\int_{B(0,c_n\sqrt{N+1})}|P(x)|^2e^{-|x|^2}dx \leq \frac{4}{3}\|P\|_{L^2(B(0,c_n\sqrt{N+1}))}^2
\end{equation}
and
\begin{multline}\label{we4x}
\|P\|_{L^2(\omega \cap B(0,c_n\sqrt{N+1}))}^2=\int_{\omega \cap B(0,c_n\sqrt{N+1})}|P(x)|^2e^{-|x|^2}e^{|x|^2}dx \\ \leq e^{c_n^2(N+1)}\|f\|_{L^2(\omega \cap B(0,c_n\sqrt{N+1}))}^2.
\end{multline}
We deduce from Lemma~\ref{remez-1} and (\ref{dfx}) that 
\begin{multline}\label{we4xx}
\|P\|_{L^2(B(0,c_n\sqrt{N+1}))}^2 \\
\leq \frac{2^{4N+2}}{3}\frac{4|B(0,c_n\sqrt{N+1})|}{|\omega \cap B(0,c_n\sqrt{N+1})|}\Big[F\Big(\frac{|\omega \cap B(0,c_n\sqrt{N+1})|}{4|B(0,c_n\sqrt{N+1})|}\Big)\Big]^{2N}
\|P\|_{L^2(\omega \cap B(0,c_n\sqrt{N+1}))}^2,
\end{multline}
with $F$ the decreasing function
$$\forall 0<t\leq 1, \quad F(t)=\frac{1+(1-t)^{\frac{1}{n}}}{1-(1-t)^{\frac{1}{n}}} \geq 1.$$
By using that $F$ is a decreasing function, it follows from (\ref{dfx}) and (\ref{we4xx}) that 
\begin{equation}\label{we5x}
\|P\|_{L^2(B(0,c_n\sqrt{N+1}))}^2 
\leq \frac{2^{4N+4}}{3\delta}\Big[F\Big(\frac{\delta}{4}\Big)\Big]^{2N}
\|P\|_{L^2(\omega \cap B(0,c_n\sqrt{N+1}))}^2.
\end{equation}
Putting together (\ref{we3x}), (\ref{we4x}) and (\ref{we5x}), we deduce that there exists a positive constant $C=C(\delta,n)>0$ such that for all $N \in \nn$ with $c_n\sqrt{N+1} \geq R_0$ and for all $f \in \mathcal E_{N}$,
\begin{equation}\label{we6x}
\|f\|_{L^2(\rr^n)}^2 \leq \frac{2^{4N+6}}{9\delta}\Big[F\Big(\frac{\delta}{4}\Big)\Big]^{2N}
e^{c_n^2(N+1)}\|f\|_{L^2(\omega \cap B(0,c_n\sqrt{N+1}))}^2\leq C^2e^{2CN}\|f\|_{L^2(\omega)}^2.
\end{equation}
The two estimates (\ref{spec1x}) and (\ref{we6x}) allow to prove the assertion $(ii)$ in Theorem~\ref{th1}.

\subsubsection{Case when the control subset is a thick set}\label{thickproof}

Let $\omega$ be a measurable subset of~$\rr^n$. We assume that $\omega$ is $\gamma$-thick at scale $L>0$, 
\begin{equation}\label{thick}
\exists 0<\gamma \leq 1, \exists L>0, \forall x \in \rr^n, \quad |\omega \cap (x+[0,L]^n)| \geq \gamma L^n.
\end{equation}
The following proof is an adaptation of the proof of the sharpened version of the Logvinenko-Sereda theorem given by Kovrijkine in~\cite{thesis,Kovrijkine}.

\medskip

\noindent
\textsc{Step 1. Bad and good cubes.} Let $N \in \nn$ be a non-negative integer and $f \in \mathcal E_{N} \setminus \{0\}$. For each multi-index $\alpha=(\alpha_1,...,\alpha_n) \in (L\mathbb{Z})^n$, let
$$Q(\alpha)=\Big\{x=(x_1,...,x_n) \in \rr^n :\ \forall 1 \leq j \leq n, \ |x_j-\alpha_j| <\frac{L}{2}\Big\}.$$
Notice that
$$\forall \alpha, \beta \in (L\mathbb{Z})^n, \ \alpha \neq \beta, \quad Q(\alpha) \cap Q(\beta)=\emptyset, \qquad \rr^n=\bigcup_{\alpha \in (L\mathbb{Z})^n}\overline{Q(\alpha)},$$
where $\overline{Q(\alpha)}$ denotes the closure of $Q(\alpha)$. It follows that
$$\|f\|_{L^2(\rr^n)}^2=\sum_{\alpha \in (L\mathbb{Z})^n}\int_{Q(\alpha)}|f(x)|^2dx.$$
Let $\delta>0$ be a positive constant to be chosen later on.
We divide the family of cubes $(Q(\alpha))_{\alpha \in (L\mathbb{Z})^n}$ into families of good and bad cubes. A cube $Q(\alpha)$, with $\alpha \in (L\mathbb{Z})^n$, is said to be good if it satisfies for all $\beta \in \nn^n$,  
\begin{equation}\label{good}
\int_{Q(\alpha)}|\partial_x^{\beta}f(x)|^2dx\leq e^{e\delta^{-2}}\big(8\delta^2(2^n+1)\big)^{|\beta|}
(|\beta|!)^2e^{2\delta^{-1}\sqrt{N}} \int_{Q(\alpha)}|f(x)|^2dx.
\end{equation}
On the other hand, a cube $Q(\alpha)$, with $\alpha \in (L\mathbb{Z})^n$, which is not good, is said to be bad, that is,
\begin{multline}\label{bad}
\exists \beta \in \nn^n, \ |\beta|>0, \\  \int_{Q(\alpha)}|\partial_x^{\beta}f(x)|^2dx > e^{e\delta^{-2}}\big(8\delta^2(2^n+1)\big)^{|\beta|}
(|\beta|!)^2e^{2\delta^{-1}\sqrt{N}} \int_{Q(\alpha)}|f(x)|^2dx.
\end{multline}
If $Q(\alpha)$ is a bad cube, it follows from (\ref{bad}) that there exists $\beta_0 \in \nn^n$, $|\beta_0|>0$ such that 
\begin{multline}\label{gh5}
\int_{Q(\alpha)}|f(x)|^2dx \leq  
\frac{e^{-e\delta^{-2}}}{\big(8\delta^2(2^{n}+1)\big)^{|\beta_0|}(|\beta_0|!)^2e^{2\delta^{-1}\sqrt{N}}}\int_{Q(\alpha)}|\partial_x^{\beta_0}f(x)|^2dx \\
\leq \sum_{\beta \in \nn^n,  |\beta|>0} \frac{e^{-e\delta^{-2}}}{\big(8\delta^2(2^{n}+1)\big)^{|\beta|}
(|\beta|!)^2e^{2\delta^{-1}\sqrt{N}}}\int_{Q(\alpha)}|\partial_x^{\beta}f(x)|^2dx.
\end{multline}
By summing over all the bad cubes, we deduce from (\ref{gh5}) and the Fubini-Tonelli theorem that 
\begin{align}\label{gh6}
& \ \int_{\bigcup_{\textrm{bad cubes}} Q(\alpha)}|f(x)|^2dx=\sum_{\textrm{bad cubes}}\int_{Q(\alpha)}|f(x)|^2dx \\ \notag
\leq & \ \sum_{\beta \in \nn^n,  |\beta|>0} \frac{e^{-e\delta^{-2}}}{\big(8\delta^2(2^{n}+1)\big)^{|\beta|}
(|\beta|!)^2e^{2\delta^{-1}\sqrt{N}}}\int_{\bigcup_{\textrm{bad cubes}}  Q(\alpha)}|\partial_x^{\beta}f(x)|^2dx\\  \notag
\leq & \ \sum_{\beta \in \nn^n,  |\beta|>0} \frac{e^{-e\delta^{-2}}}{\big(8\delta^2(2^{n}+1)\big)^{|\beta|}
(|\beta|!)^2e^{2\delta^{-1}\sqrt{N}}}\int_{\rr^n}|\partial_x^{\beta}f(x)|^2dx.
\end{align}
By using that the number of solutions to the equation $\beta_1+...+\beta_{n}=k$, with $k \geq 0$, $n \geq 1$ and unknown $\beta=(\beta_1,...,\beta_n) \in \nn^{n}$, is given by $\binom{k+n-1}{k}$, we obtain from the Bernstein type estimates in Proposition~\ref{prop1} (formula $(i)$) and (\ref{gh6}) that 
\begin{align}\label{gh6y}
\int_{\bigcup_{\textrm{bad cubes}} Q(\alpha)}|f(x)|^2dx
\leq &\ \Big(\sum_{\beta \in \nn^n,  |\beta|>0} \frac{1}{\big(2(2^{n}+1)\big)^{|\beta|}}\Big)
\|f\|_{L^2(\rr^n)}^2\\ \notag
= &\ \Big(\sum_{k=1}^{+\infty}\binom{k+n-1}{k} \frac{1}{2^k(2^{n}+1)^{k}}\Big)
\|f\|_{L^2(\rr^n)}^2\\ \notag
\leq & \ 2^{n-1}\Big(\sum_{k=1}^{+\infty}\frac{1}{(2^{n}+1)^{k}}\Big)
\|f\|_{L^2(\rr^n)}^2=\frac{1}{2}\|f\|_{L^2(\rr^n)}^2,
\end{align}
since 
\begin{equation}\label{gh45}
\binom{k+n-1}{k} \leq \sum_{j=0}^{k+n-1}\binom{k+n-1}{j}=2^{k+n-1}.
\end{equation}
By writing
$$\|f\|_{L^2(\rr^n)}^2=\int_{\bigcup_{\textrm{good cubes}} Q(\alpha)}|f(x)|^2dx+ \int_{\bigcup_{\textrm{bad cubes}} Q(\alpha)}|f(x)|^2dx,$$
it follows from (\ref{gh6y}) that 
\begin{equation}\label{gh7}
\|f\|_{L^2(\rr^n)}^2 \leq 2 \int_{\bigcup_{\textrm{good cubes}} Q(\alpha)}|f(x)|^2dx.
\end{equation}

\medskip

\noindent
\textsc{Step 2. Properties on good cubes.}
As any cube $Q(\alpha)$ satisfies the cone condition, the Sobolev embedding 
$$W^{n,2}(Q(\alpha)) \xhookrightarrow{} L^{\infty}(Q(\alpha)),$$
see e.g.~\cite[Theorem~4.12]{adams}, implies that there exists a universal positive constant $C_n>0$ depending only on the dimension $n \geq 1$ such that 
\begin{equation}\label{sobolev}
\forall u \in W^{n,2}(Q(\alpha)), \quad 
\|u\|_{L^{\infty}(Q(\alpha))} \leq C_n \|u\|_{W^{n,2}(Q(\alpha))}.
\end{equation}
By translation invariance of the Lebesgue measure, notice in particular that the constant $C_n$ does not depend on the parameter $\alpha \in (L\mathbb{Z})^n$. Let $Q(\alpha)$ be a good cube. We deduce from (\ref{good}) and (\ref{sobolev}) that for all $\beta \in \nn^n$,
\begin{align}\label{gh30}
& \ \|\partial_x^{\beta}f\|_{L^{\infty}(Q(\alpha))} 
\leq C_n\Big(\sum_{\tilde{\beta} \in \nn^n, |\tilde{\beta}| \leq n}\|\partial_x^{\beta+\tilde{\beta}}f\|^2_{L^{2}(Q(\alpha))}\Big)^{\frac{1}{2}}\\
\notag
\leq & \ C_ne^{\frac{e \delta^{-2}}{2}}e^{\delta^{-1}\sqrt{N}}\Big(\sum_{\tilde{\beta} \in \nn^n, |\tilde{\beta}| \leq n}
\big(8\delta^2(2^{n}+1)\big)^{|\beta|+|\tilde{\beta}|}
\big((|\beta|+|\tilde{\beta}|)!\big)^2 \Big)^{\frac{1}{2}}\|f\|_{L^2(Q(\alpha))}\\ \notag
\leq & \ \tilde{C}_n(\delta)\big(32\delta^2(2^{n}+1)\big)^{\frac{|\beta|}{2}}|\beta|!e^{\delta^{-1}\sqrt{N}}\|f\|_{L^2(Q(\alpha))},
\end{align}
with 
\begin{equation}\label{gh50}
\tilde{C}_n(\delta)=C_ne^{\frac{e \delta^{-2}}{2}}\Big(\sum_{\tilde{\beta} \in \nn^n, |\tilde{\beta}| \leq n}
\big(32\delta^2(2^{n}+1)\big)^{|\tilde{\beta}|}(|\tilde{\beta}|!)^2 \Big)^{\frac{1}{2}}>0,
\end{equation}
since 
$$(|\beta|+|\tilde{\beta}|)! \leq 2^{|\beta|+|\tilde{\beta}|}|\beta|!|\tilde{\beta}|!.$$
Recalling that $f$ is a finite combination of Hermite functions, we deduce from the continuity of the function $f$ and the compactness of $\overline{Q(\alpha)}$ that there exists $x_{\alpha} \in \overline{Q(\alpha)}$ such that 
\begin{equation}\label{gh8}
\|f\|_{L^{\infty}(Q(\alpha))}=|f(x_{\alpha})|.
\end{equation}
By using spherical coordinates centered at $x_{\alpha} \in \overline{Q(\alpha)}$ and the fact that the Euclidean diameter of the cube $Q(\alpha)$ is $\sqrt{n}L$, we observe that 
\begin{align}\label{gh9}
|\omega \cap Q(\alpha)|=& \ \int_0^{+\infty}\Big(\int_{\mathbb{S}^{n-1}}\un_{\omega \cap Q(\alpha)}(x_{\alpha}+r \sigma)d\sigma\Big)r^{n-1}dr\\ \notag
= & \ \int_0^{\sqrt{n}L}\Big(\int_{\mathbb{S}^{n-1}}\un_{\omega \cap Q(\alpha)}(x_{\alpha}+r \sigma)d\sigma\Big)r^{n-1}dr\\ \notag
= & \ n^{\frac{n}{2}}L^n\int_0^{1}\Big(\int_{\mathbb{S}^{n-1}}\un_{\omega \cap Q(\alpha)}(x_{\alpha}+\sqrt{n}Lr \sigma)d\sigma\Big)r^{n-1}dr,
\end{align}
where $\un_{\omega \cap Q(\alpha)}$ denotes the characteristic function of the measurable set $\omega \cap Q(\alpha)$.
By using the Fubini's theorem, we deduce from (\ref{gh9}) that 
\begin{align}\label{gh10}
|\omega \cap Q(\alpha)| \leq & \ n^{\frac{n}{2}}L^n\int_0^{1}\Big(\int_{\mathbb{S}^{n-1}}\un_{\omega \cap Q(\alpha)}(x_{\alpha}+\sqrt{n}Lr \sigma)d\sigma\Big)dr\\ \notag
= & \ n^{\frac{n}{2}}L^n\int_{\mathbb{S}^{n-1}}\Big(\int_0^{1}\un_{\omega \cap Q(\alpha)}(x_{\alpha}+\sqrt{n}Lr \sigma)dr\Big)d\sigma\\ \notag
= & \ n^{\frac{n}{2}}L^n\int_{\mathbb{S}^{n-1}}\Big(\int_0^{1}\un_{I_{\sigma}}(r)dr\Big)d\sigma=n^{\frac{n}{2}}L^n\int_{\mathbb{S}^{n-1}}|I_{\sigma}|d\sigma,
\end{align}
where
\begin{equation}\label{gh11}
I_{\sigma}=\{r \in [0,1] :\ x_{\alpha}+\sqrt{n}Lr \sigma \in \omega \cap Q(\alpha)\}.
\end{equation}
The estimate (\ref{gh10}) implies that there exists $\sigma_0 \in \mathbb{S}^{n-1}$ such that 
\begin{equation}\label{gh12}
|\omega \cap Q(\alpha)| \leq n^{\frac{n}{2}}L^n|\mathbb{S}^{n-1}| |I_{\sigma_0}|.
\end{equation}
By using the thickness property (\ref{thick}), it follows from (\ref{gh12}) that 
\begin{equation}\label{gh13}
|I_{\sigma_0}| \geq \frac{\gamma}{n^{\frac{n}{2}}|\mathbb{S}^{n-1}|}>0.
\end{equation}

\medskip

\noindent
\textsc{Step 3. Recovery of the $L^2(\rr^n)$-norm.}
We first notice that $\|f\|_{L^2(Q(\alpha))} \neq 0$, since $f$ is a non-zero entire function. We consider the entire function
\begin{equation}\label{gh13b}
\forall z \in \cc, \quad \phi(z)=L^{\frac{n}{2}}\frac{f(x_{\alpha}+\sqrt{n}Lz \sigma_0)}{\|f\|_{L^2(Q(\alpha))}}.
\end{equation}
We observe from (\ref{gh8}) that 
$$|\phi(0)|=\frac{L^{\frac{n}{2}}|f(x_{\alpha})|}{\|f\|_{L^2(Q(\alpha))}}=\frac{L^{\frac{n}{2}}\|f\|_{L^{\infty}(Q(\alpha))}}{\|f\|_{L^2(Q(\alpha))}} \geq 1.$$
Instrumental in the proof is the following lemma proved by Kovrijkine in~\cite[Lemma~1]{Kovrijkine}:

\medskip

\begin{lemma}\label{lem_kov} \emph{(\cite[Lemma~1]{Kovrijkine})}.
Let $I \subset \rr$ be an interval of length $1$ such that $0 \in I$ and $E \subset I$ be a subset of positive measure $|E|>0$. There exists a positive constant $C>1$ such that for all analytic function $\Phi$ on the open disc $D(0,5)$ such that $|\Phi(0)| \geq 1$, then 
$$\sup_{x \in I}|\Phi(x)| \leq \Big(\frac{C}{|E|}\Big)^{\frac{\ln M}{\ln 2}}\sup_{x \in E}|\Phi(x)| ,$$
with $M=\sup_{|z| \leq 4}|\Phi(z)| \geq 1$.  
\end{lemma}

\medskip

Applying Lemma~\ref{lem_kov} with $I=[0,1]$, $E=I_{\sigma_0} \subset [0,1]$ verifying $|E|=|I_{\sigma_0}|>0$ according to (\ref{gh13}), and the analytic function $\Phi=\phi$ defined in (\ref{gh13b}) satisfying $|\phi(0)| \geq 1$, we obtain that
\begin{equation}\label{gh20}
L^{\frac{n}{2}}\frac{\sup_{x \in [0,1]}|f(x_{\alpha}+\sqrt{n}Lx \sigma_0)|}{\|f\|_{L^2(Q(\alpha))}} \leq \Big(\frac{C}{|I_{\sigma_0}|}\Big)^{\frac{\ln M}{\ln 2}}L^{\frac{n}{2}}\frac{\sup_{x \in I_{\sigma_0}}|f(x_{\alpha}+\sqrt{n}Lx \sigma_0)|}{\|f\|_{L^2(Q(\alpha))}},
\end{equation}
with
\begin{equation}\label{gh22}
M=L^{\frac{n}{2}}\frac{\sup_{|z| \leq 4}|f(x_{\alpha}+\sqrt{n}Lz \sigma_0)|}{\|f\|_{L^2(Q(\alpha))}}.
\end{equation}
It follows from (\ref{gh13}) and (\ref{gh20}) that 
\begin{multline}\label{gh21}
\sup_{x \in [0,1]}|f(x_{\alpha}+\sqrt{n}Lx \sigma_0)| \leq \Big(\frac{C n^{\frac{n}{2}}|\mathbb{S}^{n-1}|}{\gamma}\Big)^{\frac{\ln M}{\ln 2}}\sup_{x \in I_{\sigma_0}}|f(x_{\alpha}+\sqrt{n}Lx \sigma_0)|\\
\leq M^{\frac{1}{\ln 2}\ln(\frac{Cn^{\frac{n}{2}}|\mathbb{S}^{n-1}|}{\gamma})}
\sup_{x \in I_{\sigma_0}}|f(x_{\alpha}+\sqrt{n}Lx \sigma_0)|.
\end{multline}
According to (\ref{gh11}), we notice that 
\begin{equation}\label{gh23}
\sup_{x \in I_{\sigma_0}}|f(x_{\alpha}+\sqrt{n}Lx \sigma_0)| \leq \|f\|_{L^{\infty}(\omega \cap Q(\alpha))}.
\end{equation}
On the other hand, we deduce from (\ref{gh8}) that 
\begin{equation}\label{gh24}
\|f\|_{L^{\infty}(Q(\alpha))}= |f(x_{\alpha})| \leq \sup_{x \in [0,1]}|f(x_{\alpha}+\sqrt{n}Lx \sigma_0)|.
\end{equation}
It follows from (\ref{gh21}), (\ref{gh23}) and (\ref{gh24}) that 
\begin{equation}\label{gh25}
\|f\|_{L^{\infty}(Q(\alpha))}\leq M^{\frac{1}{\ln 2}\ln(\frac{Cn^{\frac{n}{2}}|\mathbb{S}^{n-1}|}{\gamma})}\|f\|_{L^{\infty}(\omega \cap Q(\alpha))}.
\end{equation}
By using the analyticity of the function $f$, we observe that 
\begin{equation}\label{gh31}
\forall z \in \cc, \quad f(x_{\alpha}+\sqrt{n}Lz \sigma_0)
=\sum_{\beta \in \nn^n}\frac{(\partial_x^{\beta}f)(x_{\alpha})}{\beta!}\sigma_0^{\beta}n^{\frac{|\beta|}{2}}L^{|\beta|}z^{|\beta|}.
\end{equation}
By using that $Q(\alpha)$ is a good cube, $x_{\alpha} \in \overline{Q(\alpha)}$ and the continuity of the functions $\partial_x^{\beta}f$, we deduce from (\ref{gh30}) and (\ref{gh31}) that for all $|z| \leq 4$,
\begin{multline}\label{gh32}
|f(x_{\alpha}+\sqrt{n}Lz \sigma_0)|
\leq \sum_{\beta \in \nn^n}\frac{|(\partial_x^{\beta}f)(x_{\alpha})|}{\beta!}(4\sqrt{n}L)^{|\beta|} \\
\leq \tilde{C}_n(\delta)e^{\delta^{-1}\sqrt{N}}\Big(\sum_{\beta \in \nn^n}\frac{|\beta|!}{\beta!}\big(\delta L\sqrt{2^{9}n(2^{n}+1)}\big)^{|\beta|}\Big)\|f\|_{L^2(Q(\alpha))}.
\end{multline}
By using anew that the number of solutions to the equation $\beta_1+...+\beta_{n}=k$, with $k \geq 0$, $n \geq 1$ and unknown $\beta=(\beta_1,...,\beta_n) \in \nn^{n}$, is given by $\binom{k+n-1}{k}$, and that
$$|\beta|! \leq n^{|\beta|}\beta!,$$
see e.g. (\ref{esti3}), we notice from (\ref{gh45}) that 
\begin{multline}\label{gh49}
\sum_{\beta \in \nn^n}\frac{|\beta|!}{\beta!}\big(\delta L\sqrt{2^{9}n(2^{n}+1})\big)^{|\beta|} \leq \sum_{\beta \in \nn^n}\big(\delta L\sqrt{2^9n^3(2^{n}+1)}\big)^{|\beta|}\\
=\sum_{k=0}^{+\infty}\binom{k+n-1}{k}\big(\delta L\sqrt{2^9n^3(2^{n}+1)}\big)^{k}
\leq 2^{n-1}\sum_{k=0}^{+\infty}\big(\delta L\sqrt{2^{11}n^3(2^{n}+1)}\big)^{k}.
\end{multline}
From now on, the positive parameter $\delta>0$ is fixed and taken to be equal to 
\begin{equation}\label{gh51}
\delta=\frac{1}{\delta_nL}>0,
\end{equation}
with
$$\delta_n=2\sqrt{2^{11}n^3(2^{n}+1)}>0.$$
With this choice, it follows from (\ref{gh22}), (\ref{gh32}), (\ref{gh49}) and (\ref{gh51}) that
\begin{equation}\label{gh52}
M \leq (4L)^{\frac{n}{2}}\tilde{C}_n(\delta_n^{-1}L^{-1})e^{\delta_n L\sqrt{N}}.
\end{equation}
The positive constant $C>1$ given by Lemma~\ref{lem_kov} may be chosen such that
\begin{equation}\label{gh52s}
Cn^{\frac{n}{2}}|\mathbb{S}^{n-1}| >1.
\end{equation}  
With this choice, we deduce from (\ref{gh25}) and (\ref{gh52}) that 
\begin{equation}\label{gh53}
\|f\|_{L^{\infty}(Q(\alpha))}\leq 
\Big(\frac{Cn^{\frac{n}{2}}|\mathbb{S}^{n-1}|}{\gamma}\Big)^{\frac{\ln((4L)^{\frac{n}{2}}\tilde{C}_n(\delta_n^{-1}L^{-1}))}{\ln 2}+\frac{\delta_n}{\ln 2}L\sqrt{N}}
\|f\|_{L^{\infty}(\omega \cap Q(\alpha))}.
\end{equation}
Recalling from the thickness property (\ref{thick}) that $|\omega \cap Q(\alpha)| \geq \gamma L^n>0$ and
setting 
\begin{equation}\label{gh56}
\tilde{\omega}_{\alpha}=\Big\{x \in \omega \cap Q(\alpha) :\ |f(x)| \leq \frac{2}{|\omega \cap Q(\alpha)|}\int_{\omega \cap Q(\alpha)}|f(t)|dt\Big\},
\end{equation}
we observe that 
\begin{multline}\label{gh57}
\int_{\omega \cap Q(\alpha)}|f(x)|dx \geq \int_{(\omega \cap Q(\alpha))\setminus \tilde{\omega}_{\alpha}}|f(x)|dx \geq \frac{2|(\omega \cap Q(\alpha))\setminus \tilde{\omega}_{\alpha}|}{|\omega \cap Q(\alpha)|}\int_{\omega \cap Q(\alpha)}|f(x)|dx.
\end{multline}
By using that the integral 
$$\int_{\omega \cap Q(\alpha)}|f(x)|dx > 0,$$
since $f$ is a non-zero entire function and $|\omega \cap Q(\alpha)|>0$, we obtain that
$$|(\omega \cap Q(\alpha))\setminus \tilde{\omega}_{\alpha}| \leq \frac{1}{2}|\omega \cap Q(\alpha)|,$$
which implies that 
\begin{equation}\label{gh58}
|\tilde{\omega}_{\alpha}|=|\omega \cap Q(\alpha)|- |(\omega \cap Q(\alpha))\setminus \tilde{\omega}_{\alpha}| \geq \frac{1}{2}|\omega \cap Q(\alpha)|\geq \frac{1}{2}\gamma L^n>0,
\end{equation}
thanks anew to the thickness property (\ref{thick}).
By using again spherical coordinates as in (\ref{gh9}) and (\ref{gh10}), we observe that 
\begin{multline}\label{gh9a}
|\tilde{\omega}_{\alpha}|=|\tilde{\omega}_{\alpha} \cap Q(\alpha)|\\
=n^{\frac{n}{2}}L^n\int_0^{1}\Big(\int_{\mathbb{S}^{n-1}}\un_{\tilde{\omega}_{\alpha}  \cap Q(\alpha)}(x_{\alpha}+\sqrt{n}Lr \sigma)d\sigma\Big)r^{n-1}dr  \leq n^{\frac{n}{2}}L^n\int_{\mathbb{S}^{n-1}}|\tilde{I}_{\sigma}|d\sigma,
\end{multline}
where
\begin{equation}\label{gh11a}
\tilde{I}_{\sigma}=\{r \in [0,1] :\ x_{\alpha}+\sqrt{n}Lr \sigma \in \tilde{\omega}_{\alpha} \cap Q(\alpha)\}.
\end{equation}
As in (\ref{gh12}), the estimate (\ref{gh9a}) implies that there exists $\sigma_0 \in \mathbb{S}^{n-1}$ such that 
\begin{equation}\label{gh12a}
|\tilde{\omega}_{\alpha}| \leq n^{\frac{n}{2}}L^n|\mathbb{S}^{n-1}| |\tilde{I}_{\sigma_0}|.
\end{equation}
We deduce from (\ref{gh58}) and (\ref{gh12a}) that 
\begin{equation}\label{gh13a}
|\tilde{I}_{\sigma_0}| \geq \frac{\gamma}{2n^{\frac{n}{2}}|\mathbb{S}^{n-1}|}>0.
\end{equation}
Applying anew Lemma~\ref{lem_kov} with $I=[0,1]$, $E=\tilde{I}_{\sigma_0} \subset [0,1]$ verifying $|E|=|\tilde{I}_{\sigma_0}|>0$, and the analytic function $\Phi=\phi$ defined in (\ref{gh13b}) satisfying $|\phi(0)| \geq 1$, we obtain that
\begin{equation}\label{gh20a}
L^{\frac{n}{2}}\frac{\sup_{x \in [0,1]}|f(x_{\alpha}+\sqrt{n}Lx \sigma_0)|}{\|f\|_{L^2(Q(\alpha))}} \leq \Big(\frac{C}{|\tilde{I}_{\sigma_0}|}\Big)^{\frac{\ln M}{\ln 2}}L^{\frac{n}{2}}\frac{\sup_{x \in \tilde{I}_{\sigma_0}}|f(x_{\alpha}+\sqrt{n}Lx \sigma_0)|}{\|f\|_{L^2(Q(\alpha))}},
\end{equation}
where $M$ denotes the constant defined in (\ref{gh22}).
It follows from (\ref{gh13a}) and (\ref{gh20a}) that 
\begin{multline}\label{gh21a}
\sup_{x \in [0,1]}|f(x_{\alpha}+\sqrt{n}Lx \sigma_0)| \leq \Big(\frac{2C n^{\frac{n}{2}}|\mathbb{S}^{n-1}|}{\gamma}\Big)^{\frac{\ln M}{\ln 2}}\sup_{x \in \tilde{I}_{\sigma_0}}|f(x_{\alpha}+\sqrt{n}Lx \sigma_0)|\\
\leq M^{\frac{1}{\ln 2}\ln(\frac{2Cn^{\frac{n}{2}}|\mathbb{S}^{n-1}|}{\gamma})}
\sup_{x \in \tilde{I}_{\sigma_0}}|f(x_{\alpha}+\sqrt{n}Lx \sigma_0)|.
\end{multline}
According to (\ref{gh11a}), we notice that 
\begin{equation}\label{gh23a}
\sup_{x \in \tilde{I}_{\sigma_0}}|f(x_{\alpha}+\sqrt{n}Lx \sigma_0)| \leq \|f\|_{L^{\infty}(\tilde{\omega}_{\alpha} \cap Q(\alpha))}.
\end{equation}
It follows from (\ref{gh24}), (\ref{gh21a}) and (\ref{gh23a}) that 
\begin{equation}\label{gh25a}
\|f\|_{L^{\infty}(Q(\alpha))}\leq M^{\frac{1}{\ln 2}\ln(\frac{2Cn^{\frac{n}{2}}|\mathbb{S}^{n-1}|}{\gamma})}\|f\|_{L^{\infty}(\tilde{\omega}_{\alpha} \cap Q(\alpha))}.
\end{equation}
On the other hand, it follows from (\ref{gh56}) that
\begin{equation}\label{gh60}
\|f\|_{L^{\infty}(\tilde{\omega}_{\alpha} \cap Q(\alpha))} \leq \frac{2}{|\omega \cap Q(\alpha)|}\int_{\omega \cap Q(\alpha)}|f(x)|dx.
\end{equation}
We deduce from (\ref{gh25a}), (\ref{gh60}) and the Cauchy-Schwarz inequality that 
\begin{align}\label{gh61}
 \|f\|_{L^{2}(Q(\alpha))} \leq & \ L^{\frac{n}{2}} \|f\|_{L^{\infty}(Q(\alpha))}\\ \notag
\leq & \ \frac{2L^{\frac{n}{2}}}{|\omega \cap Q(\alpha)|} M^{\frac{1}{\ln 2}\ln(\frac{2Cn^{\frac{n}{2}}|\mathbb{S}^{n-1}|}{\gamma})}
\int_{\omega \cap Q(\alpha)}|f(x)|dx\\ \notag
\leq & \  \frac{2L^{\frac{n}{2}}}{|\omega \cap Q(\alpha)|^{\frac{1}{2}}} M^{\frac{1}{\ln 2}\ln(\frac{2Cn^{\frac{n}{2}}|\mathbb{S}^{n-1}|}{\gamma})}
\|f\|_{L^2(\omega \cap Q(\alpha))}.
\end{align}
By using the thickness property (\ref{thick}), it follows from (\ref{gh52}), (\ref{gh52s}) and (\ref{gh61}) that
\begin{multline}\label{gh62}
\|f\|_{L^{2}(Q(\alpha))}^2 \leq  \frac{4}{\gamma} M^{\frac{2}{\ln 2}\ln(\frac{2Cn^{\frac{n}{2}}|\mathbb{S}^{n-1}|}{\gamma})}
\|f\|_{L^2(\omega \cap Q(\alpha))}^2\\
\leq  \frac{4}{\gamma} \Big((4L)^{\frac{n}{2}}\tilde{C}_n(\delta_n^{-1}L^{-1})e^{\delta_n L\sqrt{N}}\Big)^{\frac{2}{\ln 2}\ln(\frac{2Cn^{\frac{n}{2}}|\mathbb{S}^{n-1}|}{\gamma})}
\|f\|_{L^2(\omega \cap Q(\alpha))}^2.
\end{multline}
With 
\begin{equation}\label{gh54}
\kappa_n(L,\gamma)=\frac{2^{\frac{3}{2}}}{\gamma^{\frac{1}{2}}}
\Big(\frac{2Cn^{\frac{n}{2}}|\mathbb{S}^{n-1}|}{\gamma}\Big)^{\frac{\ln((4L)^{\frac{n}{2}}\tilde{C}_n(\delta_n^{-1}L^{-1}))}{\ln 2}}>0,
\end{equation}
we deduce from (\ref{gh62}) that there exists a positive universal constant $\tilde{\kappa}_n>0$ such that for any good cube $Q(\alpha)$,
\begin{equation}\label{gh55}
\|f\|_{L^{2}(Q(\alpha))}^2\leq \frac{1}{2}\kappa_n(L,\gamma)^2\Big(\frac{\tilde{\kappa}_n}{\gamma}\Big)^{2\tilde{\kappa}_n L\sqrt{N}}
\|f\|_{L^{2}(\omega \cap Q(\alpha))}^2.
\end{equation}
It follows from (\ref{gh7}) and (\ref{gh55}) that 
\begin{multline*}
\|f\|_{L^2(\rr^n)}^2 \leq 2 \int_{\bigcup_{\textrm{good cubes}} Q(\alpha)}|f(x)|^2dx=2 \sum_{\textrm{good cubes}}\|f\|_{L^{2}(Q(\alpha))}^2\\
\leq \kappa_n(L,\gamma)^2\Big(\frac{\tilde{\kappa}_n}{\gamma}\Big)^{2\tilde{\kappa}_n L\sqrt{N}} \sum_{\textrm{good cubes}}
\|f\|_{L^{2}(\omega \cap Q(\alpha))}^2
\end{multline*}
and 
\begin{multline*}
\|f\|_{L^2(\rr^n)}^2 
\leq \kappa_n(L,\gamma)^2\Big(\frac{\tilde{\kappa}_n}{\gamma}\Big)^{2\tilde{\kappa}_n L\sqrt{N}} \int_{\omega \cap (\bigcup_{\textrm{good cubes}} Q(\alpha))}|f(x)|^2dx \\
\leq \kappa_n(L,\gamma)^2\Big(\frac{\tilde{\kappa}_n}{\gamma}\Big)^{2\tilde{\kappa}_n L\sqrt{N}}\|f\|_{L^2(\omega)}^2.
\end{multline*}
This ends the proof of assertion $(iii)$ in Theorem~\ref{th1}.

\section{Applications to the null-controllability of quadratic equations}\label{quadratic1} 
This section presents a result of null-controllability for parabolic equations associated to a general class of hypoelliptic non-selfadjoint accretive quadratic operators from any thick set $\omega$ of $\rr^n$ in any positive time $T>0$. We begin by recalling few facts about quadratic operators.

\subsection{Miscellaneous facts about quadratic differential operators}\label{quadratic}

Quadratic operators are pseudodifferential operators defined in the Weyl quantization
\begin{equation}\label{3}
q^w(x,D_x) f(x) =\frac{1}{(2\pi)^n}\int_{\rr^{2n}}{e^{i(x-y) \cdot \xi}q\Big(\frac{x+y}{2},\xi\Big)f(y)dyd\xi}, 
\end{equation}
by symbols $q(x,\xi)$, with $(x,\xi) \in \rr^{n} \times \rr^n$, $n \geq 1$, which are complex-valued quadratic forms 
\begin{eqnarray*}
q : \rr_x^n \times \rr_{\xi}^n &\rightarrow& \cc\\
 (x,\xi) & \mapsto & q(x,\xi).
\end{eqnarray*}
These operators are non-selfadjoint differential operators in general; with simple and fully explicit expression since the Weyl quantization of the quadratic symbol $x^{\alpha} \xi^{\beta}$, with $(\alpha,\beta) \in \nn^{2n}$, $|\alpha+\beta| = 2$, is the differential operator
$$\frac{x^{\alpha}D_x^{\beta}+D_x^{\beta} x^{\alpha}}{2}, \quad D_x=i^{-1}\partial_x.$$
Let $q^w(x,D_x)$ be a quadratic operator defined by the Weyl quantization (\ref{3}) of a complex-valued quadratic form $q$ on the phase space $\rr^{2n}$.
The maximal closed realization of the quadratic operator $q^w(x,D_x)$ on $L^2(\rr^n)$, that is, the operator equipped with the domain
\begin{equation}\label{dom1}
D(q^w)=\big\{f \in L^2(\rr^n) : \ q^w(x,D_x)f \in L^2(\rr^n)\big\},
\end{equation}
where $q^w(x,D_x)f$ is defined in the distribution sense, is known to coincide with the graph closure of its restriction to the Schwartz space~\cite[pp.~425-426]{mehler},
$$q^w(x,D_x) : \mathscr{S}(\rr^n) \rightarrow \mathscr{S}(\rr^n).$$
Let $q : \rr_x^n \times \rr_{\xi}^n \rightarrow \mathbb{C}$
be a quadratic form defined on the phase space and write $q(\cdot,\cdot)$ for its associated polarized form. Classically, one
associates to $q$ a matrix $F \in M_{2n}(\CC)$ called its Hamilton map, or its fundamental matrix.
With $\sigma$ standing for the standard symplectic form
\begin{equation}\label{11}
\sigma((x,\xi),(y,\eta))=\langle \xi, y \rangle -\langle x, \eta\rangle=\sum_{j=1}^n(\xi_j y_j-x_j \eta_j),
\end{equation}
with $x=(x_1,...,x_n)$, $y=(y_1,....,y_n)$, $\xi=(\xi_1,...,\xi_n)$, $\eta=(\eta_1,...,\eta_n) \in \cc^n$, the Hamilton map $F$ is defined as the unique matrix satisfying the identity
\begin{equation}\label{10}
\forall  (x,\xi) \in \R^{2n},\forall (y,\eta) \in \R^{2n}, \quad q((x,\xi),(y,\eta))=\sigma((x,\xi),F(y,\eta)). 
\end{equation}
We observe from the definition that 
$$F=\frac{1}{2}\left(\begin{array}{cc}
\nabla_{\xi}\nabla_x q & \nabla_{\xi}^2q  \\
-\nabla_x^2q & -\nabla_{x}\nabla_{\xi} q 
\end{array} \right),$$
where the matrices $\nabla_x^2q=(a_{i,j})_{1 \leq i,j \leq n}$,  $\nabla_{\xi}^2q=(b_{i,j})_{1 \leq i,j \leq n}$, $\nabla_{\xi}\nabla_x q =(c_{i,j})_{1 \leq i,j \leq n}$,
$\nabla_{x}\nabla_{\xi} q=(d_{i,j})_{1 \leq i,j \leq n}$ are defined by the entries
$$a_{i,j}=\partial_{x_i,x_j}^2 q, \quad b_{i,j}=\partial_{\xi_i,\xi_j}^2q, \quad c_{i,j}=\partial_{\xi_i,x_j}^2q, \quad d_{i,j}=\partial_{x_i,\xi_j}^2q.$$
The notion of singular space was introduced in~\cite{kps2} by Hitrik and the third author by pointing out the existence of a particular vector subspace in the phase space $S \subset \rr^{2n}$, which is intrinsically associated with a given quadratic symbol~$q$. This vector subspace  
is defined as the following finite intersection of kernels
\begin{equation}\label{h1bis}
S=\Big( \bigcap_{j=0}^{2n-1}\textrm{Ker}
\big[\textrm{Re }F(\textrm{Im }F)^j \big]\Big)\cap \rr^{2n},
\end{equation}
where $\textrm{Re }F$ and $\textrm{Im }F$ stand respectively for the real and imaginary parts of the Hamilton map $F$ associated with the quadratic symbol $q$,
$$\textrm{Re }F=\frac{1}{2}(F+\overline{F}), \quad \textrm{Im }F=\frac{1}{2i}(F-\overline{F}).$$
As pointed out in \cite{kps2,kps21}, the notion of singular space plays a basic role in the understanding of the spectral and hypoelliptic properties of the (possibly) non-elliptic quadratic operator $q^w(x,D_x)$, as well as the spectral and pseudospectral properties of certain classes of degenerate doubly characteristic pseudodifferential operators~\cite{kps3}. In particular, the work~\cite[Theorem~1.2.2]{kps2} gives a complete description for the spectrum of any non-elliptic quadratic operator $q^w(x,D_x)$ whose Weyl symbol $q$ has a non-negative real part $\textrm{Re }q \geq 0$, and satisfies a condition of partial ellipticity along its singular space~$S$,
\begin{equation}\label{sm2}
(x,\xi) \in S, \quad q(x,\xi)=0 \Rightarrow (x,\xi)=0. 
\end{equation}
Under these assumptions, the spectrum of the quadratic operator $q^w(x,D_x)$ is shown to be composed of a countable number of eigenvalues with finite algebraic multiplicities. The structure of this spectrum is similar to the one known for elliptic quadratic operators~\cite{sjostrand}. This condition of partial ellipticity is generally weaker than the condition of ellipticity, $S \subsetneq \rr^{2n}$, and allows one to deal with more degenerate situations.

An important class of quadratic operators satisfying condition (\ref{sm2}) are those with zero singular spaces $S=\{0\}$. In this case, the condition of partial ellipticity trivially holds.
More specifically, these quadratic operators have been shown in \cite[Theorem~1.2.1]{kps21} to be hypoelliptic and to enjoy global subelliptic estimates of the type
\begin{multline}\label{lol1}
\exists C>0, \forall f \in \mathscr{S}(\rr^n), \\ \|\langle(x,D_x)\rangle^{2(1-\delta)} f\|_{L^2(\rr^n)} \leq C(\|q^w(x,D_x) f\|_{L^2(\rr^n)}+\|f\|_{L^2(\rr^n)}),
\end{multline}
where $\langle(x,D_x)\rangle^{2}=1+|x|^2+|D_x|^2$, with a sharp loss of derivatives $0 \leq \delta<1$ with respect to the elliptic case (case $\delta=0$), which can be explicitly derived from the structure of the singular space.

When the quadratic symbol $q$ has a non-negative real part $\textrm{Re }q \geq 0$, the singular space can be also defined in an equivalent way as the subspace in the phase space where all the Poisson brackets 
$$H_{\textrm{Im}q}^k \textrm{Re }q=\left(\frac{\partial \textrm{Im }q}{\partial\xi}\cdot \frac{\partial}{\partial x}-\frac{\partial \textrm{Im }q}{\partial x}\cdot \frac{\partial}{\partial \xi}\right)^k \textrm{Re } q, \quad k \geq 0,$$ 
are vanishing
$$S=\big\{X=(x,\xi) \in \rr^{2n} : \ (H_{\textrm{Im}q}^k \textrm{Re } q)(X)=0,\ k \geq 0\big\}.$$
This dynamical definition shows that the singular space corresponds exactly to the set of points $X \in \rr^{2n}$, where the real part of the symbol $\textrm{Re }q$ under the flow of the Hamilton vector field $H_{\textrm{Im}q}$ associated with its imaginary part
\begin{equation}\label{evg5}
t \mapsto \textrm{Re }q(e^{tH_{\textrm{Im}q}}X),
\end{equation}
vanishes to infinite order at $t=0$. This is also equivalent to the fact that the function \eqref{evg5} is identically zero on~$\rr$.

\subsection{Null-controllability of hypoelliptic quadratic equations}
We study the class of quadratic operators whose Weyl symbols have non-negative real parts $\textrm{Re }q \geq 0$, and zero singular spaces $S=\{0\}$. 
According to the above description of the singular space, these quadratic operators are exactly those whose Weyl symbols have a non-negative real part $\textrm{Re }q \geq 0$, 
becoming positive definite 
\begin{equation}\label{evg41}
\forall\,T>0, \quad \langle \textrm{Re }q\rangle_T(X)=\frac{1}{2T}\int_{-T}^T{(\textrm{Re }q)(e^{tH_{\textrm{Im}q}}X)dt} \gg 0, 
\end{equation}
after averaging by the linear flow of the Hamilton vector field associated with its imaginary part. The above notation $a \gg 0$ denotes that the quadratic form $a$ is positive definite.
These quadratic operators are also known~\cite[Theorem~1.2.1]{kps2} to generate strongly continuous contraction semigroups $(e^{-tq^w})_{t \geq 0}$ on $L^2(\rr^n)$, which are smoothing in the Schwartz space for any positive time
$$\forall t>0, \forall f \in L^2(\rr^n), \quad e^{-t q^w}f \in \mathscr{S}(\rr^n).$$
In the work~\cite[Theorem~1.2]{HPSVII}, these regularizing properties were sharpened and 
these contraction semigroups were shown to be actually smoothing for any positive time in the Gelfand-Shilov space $S_{1/2}^{1/2}(\rr^n)$: $\exists C>0$, $\exists t_0 > 0$, $\forall f \in L^2(\rr^n)$, $\forall \alpha, \beta \in \nn^n$, $\forall 0<t \leq t_0$,
\begin{equation}\label{eq1.10}
\|x^{\alpha}\partial_x^{\beta}(e^{-tq^w}f)\|_{L^{\infty}(\rr^n)} \leq \frac{C^{1+|\alpha|+|\beta|}}{t^{\frac{2k_0+1}{2}(|\alpha|+|\beta|+2n+s)}}(\alpha!)^{1/2}(\beta!)^{1/2}\|f\|_{L^2(\rr^n)},
\end{equation}
where $s$ is a fixed integer verifying $s > n/2$, and where $0 \leq k_0 \leq 2n-1$ is the smallest integer satisfying 
\begin{equation}\label{h1bis2}
\Big( \bigcap_{j=0}^{k_0}\textrm{Ker}\big[\textrm{Re }F(\textrm{Im }F)^j \big]\Big)\cap \rr^{2n}=\{0\}.
\end{equation}
Thanks to this Gelfand-Shilov smoothing effect (\ref{eq1.10}), the first and third authors have established in~\cite[Proposition~4.1]{KK1} that, for any quadratic form $q : \rr_{x,\xi}^{2n} \rightarrow \cc$ with a non-negative real part $\textrm{Re }q \geq 0$ and a zero singular space $S=\{0\}$, the dissipation estimate (\ref{eq6}) holds 
with $0 \leq k_0 \leq 2n-1$ being the smallest integer satisfying (\ref{h1bis2}). Let $\omega \subset \rr^n$ be a measurable $\gamma$-thick set at scale $L>0$. 
We can then deduce from Theorem~\ref{Meta_thm_AdaptedLRmethod} with the following choices of parameters:
\begin{itemize}
\item[$(i)$] $\Omega=\mathbb{R}^n$,
\item[$(ii)$] $A=-q^w(x,D_x)$,
\item[$(iii)$] $a=\frac{1}{2}$, $b=1$, 
\item[$(iv)$] $t_0>0$ as in (\ref{eq6}) and (\ref{eq5}),
\item[$(v)$] $m=2k_0+1$, where $k_0$ is defined in (\ref{h1bis2}),
\item[$(vi)$] any constant $c_1>0$ satisfying 
$$\forall k \geq 1,  \quad C\Big(\frac{\kappa}{\gamma}\Big)^{\kappa L\sqrt{k}} \leq e^{c_1 \sqrt{k}},$$ where the positive constants $C=C(L,\gamma,n)>0$ and $\kappa=\kappa(n)>0$ are defined in Theorem~\ref{th1} (formula $(iii)$),
\item[$(vii)$] $c_2= \frac{1}{C_0}>0$, where $C_0>1$ is defined in (\ref{eq6}) and (\ref{eq5}),
\end{itemize}
the following observability estimate in any positive time
$$\exists C>1, \forall T>0, \forall f \in L^2(\rr^n), \quad \|e^{-Tq^w} f \|_{L^2(\rr^n)}^2 \leq C\exp\Big(\frac{C}{T^{2k_0+1}}\Big) \int_0^T \|e^{-tq^w} f \|_{L^2(\omega)}^2 dt.$$
We therefore obtain the following result of null-controllability:

\medskip

\begin{theorem}\label{th2}
Let $q : \rr_{x}^{n} \times \rr_{\xi}^n \rightarrow \cc$ be a complex-valued quadratic form with a non negative real part $\emph{\textrm{Re }}q \geq 0$, and a zero singular space $S=\{0\}$. If  $\omega$ is a measurable thick subset of $\mathbb{R}^n$, then the parabolic equation 
$$\left\lbrace \begin{array}{ll}
\partial_tf(t,x) + q^w(x,D_x)f(t,x)=u(t,x)\un_{\omega}(x)\,, \quad &  x \in \mathbb{R}^n, \ t>0,\\
f|_{t=0}=f_0 \in L^2(\rr^n),                                       &  
\end{array}\right.$$
with $q^w(x,D_x)$ being the quadratic differential operator defined by the Weyl quantization of the symbol $q$, is null-controllable from the set $\omega$ in any positive time $T>0$.
\end{theorem}

\medskip

As in~\cite{KK1}, this result of null-controllability given by Theorem~\ref{th2} applies in particular for the parabolic equation associated to the Kramers-Fokker-Planck operator 
\begin{equation} \label{eq:KFP}
K=-\Delta_v+\frac{v^2}{4}+v\partial_x-\partial_xV(x)\partial_v, \quad (x,v) \in \rr^{2},
\end{equation}
with a quadratic potential
$$V(x)=\frac{1}{2}ax^2, \quad a \in \rr^*,$$
which is an example of accretive quadratic operator with a zero singular space $S=\{0\}$. It also applies in the very same way to hypoelliptic Ornstein-Uhlenbeck equations posed in $L^2(\rr^n,\rho(x)dx)$-spaces, or to hypoelliptic Fokker-Planck equations posed in $L^2(\rr^n,\rho(x)^{-1}dx)$-spaces with respect to (gaussian) invariant measures $\rho$. Indeed, as explained in~\cite{KK1} (Sections~5 and 6) and after conjugation by $\sqrt{\rho}$ or $\sqrt{\rho}^{-1}$, this problem of null-controllability in weighted $L^2$-spaces can be rephrased as a problem of null-controllability in the flat $L^2(\rr^n,dx)$-space for which Theorem~\ref{th2}  applies.
We refer the reader to the works~\cite{KK1,kps11} for detailed discussions of various physics models whose evolution turns out to be ruled by accretive quadratic operators with zero singular spaces and to which therefore apply the above result of null-controllability.

The notion of thickness is a sufficient geometric condition for control subsets to derive the null-controllability for a general class of evolution equations associated to hypoelliptic non-selfadjoint quadratic operators that includes the harmonic heat equation. It is therefore a natural question to figure out whether this condition turns out also to be sufficient. To the best of our knowledge, there is no known necessary and sufficient geometric condition to ensure the null-controllability of these evolution equations even in the case of the harmonic heat equation. Nevertheless, one can expect that the thickness condition is not sharp. Indeed, contrary to the heat equation, the solutions of the above evolution equations do enjoy specific decay properties at infinity and one can conjecture that control subsets do not necessarily need to be distributed as much at infinity as required by the thickness condition. This conjecture will be investigated in future works.

\section{Appendix}\label{appendix}

\subsection{Gelfand-Shilov regularity}\label{gelfand}

We refer the reader to the works~\cite{gelfand,rodino1,rodino,toft} and the references herein for extensive expositions of the Gelfand-Shilov regularity theory.
The Gelfand-Shilov spaces $S_{\nu}^{\mu}(\rr^n)$, with $\mu,\nu>0$, $\mu+\nu\geq 1$, are defined as the spaces of smooth functions $f \in C^{\infty}(\rr^n)$ satisfying the estimates
$$\exists A,C>0, \quad |\partial_x^{\alpha}f(x)| \leq C A^{|\alpha|}(\alpha !)^{\mu}e^{-\frac{1}{A}|x|^{1/\nu}}, \quad x \in \rr^n, \ \alpha \in \mathbb{N}^n,$$
or, equivalently
$$\exists A,C>0, \quad \sup_{x \in \rr^n}|x^{\beta}\partial_x^{\alpha}f(x)| \leq C A^{|\alpha|+|\beta|}(\alpha !)^{\mu}(\beta !)^{\nu}, \quad \alpha, \beta \in \mathbb{N}^n.$$
These Gelfand-Shilov spaces  $S_{\nu}^{\mu}(\rr^n)$ may also be characterized as the spaces of Schwartz functions $f \in \mathscr{S}(\rr^n)$ satisfying the estimates
$$\exists C>0, \eps>0, \quad |f(x)| \leq C e^{-\eps|x|^{1/\nu}}, \quad x \in \rr^n, \qquad |\widehat{f}(\xi)| \leq C e^{-\eps|\xi|^{1/\mu}}, \quad \xi \in \rr^n.$$
In particular, we notice that Hermite functions belong to the symmetric Gelfand-Shilov space  $S_{1/2}^{1/2}(\rr^n)$. More generally, the symmetric Gelfand-Shilov spaces $S_{\mu}^{\mu}(\rr^n)$, with $\mu \geq 1/2$, can be nicely characterized through the decomposition into the Hermite basis $(\Phi_{\alpha})_{\alpha \in \mathbb{N}^n}$, see e.g. \cite[Proposition~1.2]{toft},
\begin{multline*}
f \in S_{\mu}^{\mu}(\rr^n) \Leftrightarrow f \in L^2(\rr^n), \ \exists t_0>0, \ \big\|\big(\langle f,\Phi_{\alpha}\rangle_{L^2}\exp({t_0|\alpha|^{\frac{1}{2\mu}})}\big)_{\alpha \in \mathbb{N}^n}\big\|_{l^2(\mathbb{N}^n)}<+\infty\\
\Leftrightarrow f \in L^2(\rr^n), \ \exists t_0>0, \ \|e^{t_0\mathcal{H}^{\frac{1}{2\mu}}}f\|_{L^2(\rr^n)}<+\infty,
\end{multline*}
where $\mathcal{H}=-\Delta_x+|x|^2$ stands for the harmonic oscillator.

\subsection{Remez inequality}\label{remez0} The classical Remez inequality~\cite{remezz}, see also~\cite{erd0,erd}, is the following estimate providing a bound on the maximum of the absolute value of an arbitrary real polynomial function $P \in \rr[X]$ of degree $d$ on $[-1,1]$ by the maximum of its absolute value on any measurable subset $E \subset [-1,1]$ of positive Lebesgue measure $0<|E|<2$,
\begin{equation}\label{remez}
\sup_{[-1,1]}|P| \leq T_d\Big(\frac{4-|E|}{|E|}\Big)\sup_{E}|P|,
\end{equation}
where 
\begin{equation}\label{che1}
T_d(X)=\frac{d}{2}\sum_{k=0}^{[\frac{d}{2}]}(-1)^k\frac{(d-k-1)!}{k!(d-2k)!}2^{d-2k}X^{d-2k}=\sum_{k=0}^{[\frac{d}{2}]}\binom{d}{2k}(X^2-1)^kX^{d-2k},
\end{equation}
see e.g.~\cite[Chapter 2]{borwein}, where $[x]$ stands the integer part of $x$, denotes the $d^{\textrm{th}}$ Chebyshev polynomial function of first kind. 
We also recall from~\cite[Chapter 2]{borwein} the definition of Chebyshev polynomial functions of second kind
\begin{equation}\label{bn1}
\forall d \in \nn, \quad U_d(X)=\sum_{k=0}^{[\frac{d}{2}]}(-1)^k\binom{d-k}{k}2^{d-2k}X^{d-2k}
\end{equation}
and
\begin{equation}\label{bn12}
\forall d \in \nn^*, \quad U_{d-1}(X)=\frac{1}{d}T_d'(X).
\end{equation}
The Remez inequality was extended in the multi-dimensional case in~\cite{brud}, see also~\cite[Formula $(4.1)$]{const} and~\cite{kroo}, as follows: for all convex bodies\footnote{A compact convex subset of $\rr^n$ with non-empty interior.} $K \subset \rr^n$, measurable subsets $E \subset K$ of positive Lebesgue measure $0<|E|<|K|$ and real polynomial functions $P \in \rr[X_1,...,X_n]$ of degree~$d$, the following estimate holds
\begin{equation}\label{remez1.1}
\sup_{K}|P| \leq T_d\left(\frac{1+(1-\frac{|E|}{|K|})^{\frac{1}{n}}}{1-(1-\frac{|E|}{|K|})^{\frac{1}{n}}}\right)\sup_{E}|P|.
\end{equation}
By recalling that all the zeros of the Chebyshev polynomial functions of first and second kind are simple and contained in the set $]-1,1[$, 
we observe from (\ref{che1}) and (\ref{bn12}) that the function $T_d$ is increasing on $[1,+\infty)$ and that  
\begin{multline}\label{che2}
\forall d \in \nn, \forall x \geq 1, \quad  1=T_d(1) \leq T_d(x)=\sum_{k=0}^{[\frac{d}{2}]}\binom{d}{2k}(x-1)^k(x+1)^kx^{d-2k}\\  \leq 
\sum_{k=0}^{[\frac{d}{2}]}\binom{d}{2k}x^k(x+x)^kx^{d-2k} = \sum_{k=0}^{[\frac{d}{2}]}\binom{d}{2k}2^kx^{d}\leq (2x)^d\sum_{k=0}^{[\frac{d}{2}]}2^k \leq (4x)^d,
\end{multline}
since $\binom{d}{2k} \leq \sum_{j=0}^d\binom{d}{j} = 2^d$. By using that 
$$\sup_{K}|Q| \leq \sup_{K}|\textrm{Re }Q|+\sup_{K}|\textrm{Im }Q|\quad \textrm{and} \quad \sup_{E}|\textrm{Re }Q|+\sup_{E}|\textrm{Im }Q| \leq 2\sup_{E}|Q|,$$
we deduce from (\ref{remez1.1}) and (\ref{che2}) that for all convex bodies $K \subset \rr^n$, measurable subsets $E \subset K$ of positive Lebesgue measure $0<|E|<|K|$, and complex polynomial functions $Q \in \cc[X_1,...,X_n]$ of degree~$d$,
\begin{equation}\label{remez1}
\sup_{K}|Q| \leq 2^{2d+1}\left(\frac{1+(1-\frac{|E|}{|K|})^{\frac{1}{n}}}{1-(1-\frac{|E|}{|K|})^{\frac{1}{n}}}\right)^d\sup_{E}|Q|.
\end{equation}
Thanks to this estimate, we can prove that the $L^2$-norm $\|\cdot\|_{L^2(\omega)}$ on any measurable subset $\omega \subset \rr^n$, with $n \geq 1$, of positive Lebesgue measure $|\omega|>0$ defines a norm on the finite dimensional vector space $\mathcal E_{N}$ defined in (\ref{jk1b}). Indeed, let $f$ be a function in $\mathcal{E}_{N}$ verifying $\|f\|_{L^2(\omega)}=0$, with $\omega \subset \rr^n$ a measurable subset of positive Lebesgue measure $|\omega|>0$. According to (\ref{defi}) and (\ref{jk1}), there exists a  
complex polynomial function $Q \in \cc[X_1,...,X_n]$ such that 
$$\forall (x_1,...,x_n) \in \rr^n, \quad f(x_1,...,x_n)=Q(x_1,...,x_n)e^{-\frac{x_1^2+...+x_n^2}{2}}.$$
The condition $\|f\|_{L^2(\omega)}=0$ first implies that $f=0$ almost everywhere in $\omega$, and therefore that $Q=0$ almost everywhere in $\omega$. We deduce from (\ref{remez1}) that the polynomial function $Q$ has to be zero on any convex body $K$ verifying $|K \cap \omega|>0$, and therefore is zero everywhere. We conclude that the $L^2$-norm $\|\cdot\|_{L^2(\omega)}$ actually defines a norm on the finite dimensional vector space $\mathcal E_{N}$.

On the other hand, the Remez inequality is a key ingredient in the proof of the following instrumental lemma needed for the proof of Theorem~\ref{th1}:

\medskip

\begin{lemma}\label{remez-1}
Let $R>0$ and $\omega \subset \rr^n$ be a measurable subset verifying $|\omega \cap B(0,R)|>0$. Then,  the following estimate holds for all complex polynomial functions $P \in \cc[X_1,...,X_n]$ of degree~$d$,
 $$\|P\|_{L^2(B(0,R))} 
\leq \frac{2^{2d+1}}{\sqrt{3}}\sqrt{\frac{4|B(0,R)|}{|\omega \cap B(0,R)|}}
\left(\frac{1+(1-\frac{|\omega \cap B(0,R)|}{4|B(0,R)|})^{\frac{1}{n}}}{1-(1-\frac{|\omega \cap B(0,R)|}{4|B(0,R)|})^{\frac{1}{n}}}\right)^d
\|P\|_{L^2(\omega \cap B(0,R))},$$
where $B(0,R)$ denotes the open Euclidean ball in $\rr^n$ centered at $0$ with radius $R>0$. 
\end{lemma}

\medskip

\begin{proof}
Let $P \in \cc[X_1,...,X_n]$ be a non-zero complex polynomial function of degree~$d$ and $R>0$. We consider the following subset 
\begin{equation}\label{bn3}
E_{\eps}=\Big\{x \in B(0,R) :\ |P(x)| \leq 2^{-2d-1}F\Big(\frac{\eps}{|B(0,R)|}\Big)^{-d}\sup_{B(0,R)}|P|\Big\},
\end{equation}
for all $0<\eps \leq B(0,R)$, and $F$ the decreasing function
\begin{equation}\label{dec}
\forall 0<t\leq 1, \quad F(t)=\frac{1+(1-t)^{\frac{1}{n}}}{1-(1-t)^{\frac{1}{n}}} \geq 1.
\end{equation}
The estimate 
$$2^{-2d-1}F\Big(\frac{\eps}{|B(0,R)|}\Big)^{-d}<1,$$
implies that $|E_{\eps}|<|B(0,R)|$.
We first check that the Lebesgue measure of this subset satisfies $|E_{\eps}| \leq \eps$. If $|E_{\eps}|>0$, it follows from (\ref{remez1}) that 
\begin{multline}\label{dec1}
0<\sup_{B(0,R)}|P| \leq 2^{2d+1}F\Big(\frac{|E_{\eps}|}{|B(0,R)|}\Big)^d\sup_{E_{\eps}}|P| \\ \leq F\Big(\frac{|E_{\eps}|}{|B(0,R)|}\Big)^dF\Big(\frac{\eps}{|B(0,R)|}\Big)^{-d}\sup_{B(0,R)}|P|.
\end{multline}
We obtain from (\ref{dec1}) that 
\begin{equation}\label{dec2}
F\Big(\frac{\eps}{|B(0,R)|}\Big) \leq F\Big(\frac{|E_{\eps}|}{|B(0,R)|}\Big). 
\end{equation}
As $F$ is a decreasing function, we deduce from (\ref{dec2}) that 
\begin{equation}\label{dec3}
\forall 0<\eps \leq B(0,R), \quad |E_{\eps}| \leq \eps. 
\end{equation}
Let $\omega \subset \rr^n$ be a measurable subset verifying $|\omega \cap B(0,R)|>0$. We consider the positive parameter
\begin{equation}\label{dec4}
0<\eps_0=\frac{1}{4}|\omega \cap B(0,R)|<|B(0,R)|.
\end{equation}
Setting
\begin{equation}\label{dec4.1}
G_{\eps_0}=\Big\{x \in B(0,R) :\ |P(x)| > 2^{-2d-1}F\Big(\frac{\eps_0}{|B(0,R)|}\Big)^{-d}\sup_{B(0,R)}|P|\Big\},
\end{equation}
we observe that 
\begin{multline}\label{dec5}
\int_{\omega \cap B(0,R)}|P(x)|^2dx \geq \int_{\omega \cap B(0,R)}\un_{G_{\eps_0}}(x)|P(x)|^2dx \\
\geq 2^{-4d-2}F\Big(\frac{\eps_0}{|B(0,R)|}\Big)^{-2d}\Big(\sup_{B(0,R)}|P|\Big)^2|\omega \cap G_{\eps_0}|.
\end{multline}
We deduce from (\ref{bn3}), (\ref{dec3}) and (\ref{dec4.1}) that 
\begin{align*}
|\omega \cap G_{\eps_0}|
=& \ |G_{\eps_0}| -\Big|\Big\{x \in B(0,R) \setminus \omega : |P(x)| > 2^{-2d-1}F\Big(\frac{\eps_0}{|B(0,R)|}\Big)^{-d}\sup_{B(0,R)}|P|\Big\}\Big|\\
\geq & \ (|B(0,R)|-|E_{\eps_0}|)-|B(0,R) \setminus \omega| \\
\geq &\  |B(0,R)|-\frac{1}{4}|\omega \cap B(0,R)|-(|B(0,R)|-|\omega \cap B(0,R)|),
\end{align*}
that is 
\begin{equation}\label{dec6}
|\omega \cap G_{\eps_0}| \geq \frac{3}{4}|\omega \cap B(0,R)|>0.
\end{equation}
It follows from (\ref{dec4}), (\ref{dec5}) and (\ref{dec6}) that 
\begin{multline}\label{dec7}
\|P\|_{L^2(B(0,R))}^2 \leq |B(0,R)|\Big(\sup_{B(0,R)}|P|\Big)^2 \\
\leq 2^{4d+2}\frac{4|B(0,R)|}{3|\omega \cap B(0,R)|}F\Big(\frac{|\omega \cap B(0,R)|}{4|B(0,R)|}\Big)^{2d} \int_{\omega \cap B(0,R)}|P(x)|^2dx.
\end{multline}
We deduce from (\ref{dec7}) that 
\begin{equation}\label{dec8}
\|P\|_{L^2(B(0,R))} 
\leq \frac{2^{2d+1}}{\sqrt{3}}\sqrt{\frac{4|B(0,R)|}{|\omega \cap B(0,R)|}}F\Big(\frac{|\omega \cap B(0,R)|}{4|B(0,R)|}\Big)^{d}\|P\|_{L^2(\omega \cap B(0,R))}.
\end{equation}
This ends the proof of Lemma~\ref{remez-1}.
\end{proof}

\end{document}